\documentclass[11 pt]{amsart}
\usepackage{amsmath, amsthm,amssymb,amsfonts}
\usepackage{enumitem}
\usepackage[mathscr]{euscript}
\usepackage[
  backend=biber,
  url=false,
  doi=false,
  isbn=false,
  eprint=false, style=alphabetic
]{biblatex}
\usepackage[pdfencoding=auto,
hidelinks
]{hyperref}
\usepackage{tabularx}
\usepackage{booktabs}
\usepackage{array}
\usepackage{xparse}    
\usepackage{bm}
\usepackage{txfonts}
\usepackage[dvipsnames, table]{xcolor}
\usepackage{xcolor}

\usepackage[ a4paper, top=2.3cm, bottom=1.6cm, left=2.0cm, right=2.0cm ]{geometry}

\makeatletter
\newcommand{\tagtext}[2]{
  \@bsphack
  \csname phantomsection\endcsname 
  \def\@currentlabel{#1}{\label{#2}}
  \@esphack
}
\makeatother

\usepackage{mathtools}
\mathtoolsset{showonlyrefs=true}
\usepackage{parskip}
\usepackage{url}
\usepackage{hyperref}
\hypersetup{
    colorlinks=true,
    linktoc=all,
    linkcolor=blue,
    citecolor=magenta,
    urlcolor=blue
}

\newtheorem{thm}{Theorem}[section]
\newtheorem{prop}[thm]{Proposition}
\newtheorem{lem}[thm]{Lemma}
\theoremstyle{plain}
\newtheorem{rmk}{Remark}
\newtheorem{cor}[thm]{Corollary}
\newcommand{\ind}{\mathbf{1}}
\newcommand{\sgn}{\operatorname{sgn}}
\newcommand{\meas}{\operatorname{meas}}
\newcommand{\lcm}{\operatorname{lcm}}

\newcommand{\eps}{\varepsilon}
\newcommand{\conj}[1]{\overline{#1}}
\newcommand{\sumn}{\sum \nolimits}

\usepackage{thm-restate}

\setlist[enumerate,1]{
    label=(\roman*),
}

\newcommand{\Z}{\mathbb{Z}}            
\newcommand{\Q}{\mathbb{Q}}            
\newcommand{\R}{\mathbb{R}}            
\newcommand{\C}{\mathbb{C}}

\newcommand{\Sp}{\operatorname{Sp}}

\newcommand{\GL}{\operatorname{GL}}    
\newcommand{\SL}{\operatorname{SL}}    

\newcommand{\z}{\mathfrak z}

\newcommand{\SLtwoZ}{\operatorname{SL}(2,\mathbb{Z})}

\newcommand{\GammaOne}[1]{\Gamma_1(#1)}

\usepackage{graphicx}
\numberwithin{equation}{section}

\title{Statistics of  fundamental Fourier coefficients of  cusp forms}

\author{Soumya Das and Animesh Sarkar}
\address{Department of Mathematics, Indian Institute of Science, Bangalore -- 560012, India.}

\email{soumya@iisc.ac.in}
\email{sanimesh@iisc.ac.in, animesh.math@outlook.com}
\subjclass[2020]{Primary: 11F30; Secondary: 11F11, 11F37, 11F46}
\keywords{Sign change, positive proportion, Fourier coefficients, cusp forms, Siegel modular forms}

\date{\today}

\begin{document}
\begin{abstract}
We study sign changes of Fourier coefficients of elliptic and Siegel modular forms in considerable generality. For any $N \ge 1$ and an arbitrary non-zero elliptic cusp form $f\in S_\kappa (\Gamma_1(N))$ ($\kappa \in \frac12 \Z$) with real Fourier coefficients, we show that a positive proportion of its non-zero Fourier coefficients change sign if $\kappa\ge 1$ is integral. If $\kappa\ge 1/2$ is half-integral, we get an almost optimal result.

As applications, we establish quantitative sign-change results for `fundamental' Fourier coefficients of Siegel cusp forms of arbitrary degree $n$. These improve the previous best known results on  the number of non-zero fundamental Fourier coefficients, even in degree $n=1$.
\end{abstract}
\maketitle

\section{Introduction}
Modular forms are central to number theory, and its Fourier coefficients are the main objects of focus. Statistical study about them, say for example their distribution, growth, congruences -- almost anything at all related to these objects hold interesting information about the modular form or spaces of modular forms. In this paper we will discuss one such statistical aspect: the frequency of changes of signs of these Fourier coefficients.

Substantial effort has been devoted to the quantitative questions about the sign changes of the Fourier coefficients of integral weight modular forms, with striking results. We only highlight a few such instances: first is the work by Kowalski-Lau-Soundararajan-Wu \cite{kow10} where it is shown that the signs of Hecke eigenvalues (if they are real) are sufficient to distinguish between newforms -- this is a surprising strengthening of the multiplicity-one result. Second, another landmark result is the work of Matom\"aki-Radziwiłł \cite{mat15}, who prove that a positive proportion of non-zero Hecke eigenvalues change sign. Furthermore, they show that asymptotically half of non-zero Fourier coefficients are positive and half are negative.
For more details, we refer the reader to these works and also to \cite{lau09} and the references therein. 
Among the plethora of other related results that are available in the literature, most of them deal with normalised Hecke eigenforms, as the multiplicative structure of the coefficients makes it possible to bring techniques from multiplicative number theory into the subject. 

The half-integral weight case is more delicate, as the Fourier coefficients lack comparable multiplicative properties. For Hecke eigenforms $f$ in the Kohnen's plus space,  the Waldspurger-Shimura correspondence relates the Fourier  coefficients of $f$ at the fundamental discriminants  to certain twisted central $L$-values. This  enables an interesting approach, thereby yielding strong statistical results. 
In particular, for half-integral weight Hecke eigenforms in the plus space of level $4$, Lester-Radziwi\l\l{} \cite{les21} unconditionally established at least about $X^{1-\eps}$ sign changes of the Fourier coefficients $a_f(n), \, 1 \le n \le X$ whereas \cite{jls} showed the same for arbitrary cusp forms of level $4N$ (in the plus space), with $N$ odd, square-free.
Under GRH, \cite{les21} showed a positive proportion of sign changes in the same set. For more details, we refer the reader to \cite{les21, jls} and references therein.

It is however not clear whether the above results might continue to hold without any Hecke eigenform or specific level assumptions. One of the purposes of this paper is to show that indeed these assumptions can be removed without any loss in the sharpness of the results. This generality  is interesting in its own right, but sometimes absolutely necessary. We explain why.

In fact, one  of the main motivations behind this work is to prove analogous results for holomorphic Siegel cusp forms of any degree $n$. To handle the sign changes of their Fourier coefficients one ultimately invokes results on elliptic modular forms discussed above -- and even if one starts with a Siegel Hecke eigenform (say of level one), still one usually never gets an elliptic eigenform (not even away from the level) to work with in the end. Moreover, this elliptic modular form has integral or half-integral weight according as $n$ is odd or even.
We will discuss about these issues later in the Introduction.
Similarly, if one wishes to study the oscillation of the ``error" in the formula for representation numbers of quadratic forms (cf. \cite{mat15}), one does not have Hecke eigenforms at one's disposal.

First let us describe the corresponding one-variable results on elliptic cusp forms which form the backbone for higher degrees.

\subsection{Results on elliptic cusp forms}
Let $f\in S_\kappa(\Gamma_1(N))$ ($\kappa \in \frac12 \Z$) be arbitrary. 
It has a Fourier  expansion 
\begin{equation}\label{expansion}
    f(\tau) = \sumn_{\substack{n\geq 1}} a_f(n) q^n.
\end{equation} 
\emph{All cusp forms in this article are assumed to have real Fourier coefficients unless otherwise specified.}

To state our results cleanly, we first define the sets and quantities which count the non-zero Fourier coefficients:
\begin{align}\label{eq:nfsx}
 \mathscr N_f(X;M) &:= \#\{X\leq n\leq MX~:~ a_f(n) \neq 0\} &= \# \{ n_1<n_2< \cdots < n_{R(X;M)} \}, \\
\mathscr N^{\flat}_f(X;M) &:=\#\{X\le n\le MX :\mu^2(n)=1,\ (n,2N)=1,\ a_f(n)\neq0\} &= \# \{ n_1<n_2< \cdots < n_{R^\flat(X;M)} \}.
\end{align}
The corresponding objects to measure changes of sign are:
\begin{align}\label{eq:sfsx}
\mathscr{S}_f(X;M) &:= \#\{1 \le i < R(X;M)~:~ a_f(n_i)a_f(n_{i+1})<0 \},\\
\mathscr S_f^{\flat}(X;M) &:=\#\{1 \le i < R^\flat(X;M):\mu^2(n_i)=1,\ (n_i,2N)=1,\ a_f(n_{i})a_f(n_{i+1})<0\},
\end{align}
where $n_i$ is strictly increasing.
We will write $\mathscr{S}_f(X)$ and $\mathscr S_f^{\flat}(X)$ when $M=2$, and indicate the dependence on $M$ only when necessary. In the paper, \textit{$M$ will be a fixed positive constant depending only on $f$.} 

Apart from its own interest, the reason for considering odd, square-free $n$ is to enable application to `fundamental' Fourier coefficients of Siegel modular forms. 
These objects are of particular significance in the theory of automorphic representations  for $\mathrm{GSp}(n)$.

To study the size of $\mathscr{S}_f(X)$ (resp. $\mathscr S_f^{\flat}(X)$), it only makes sense to compare it with $\mathscr{N}_f(X)$ (resp. $\mathscr N_f^{\flat}(X)$). 
The quantity $\mathscr{N}_f(X)$ has been widely studied in the literature and we know its size  when $\kappa \in \Z$. Indeed, from \cite[Théorème~17]{Serre81}, for $\kappa \ge 2$, 
\begin{align}\label{eq:serre-estimates}
    \mathscr{N}_f(X) \asymp \begin{cases} 
    \begin{aligned}
    &X (\log X)^{-1/2} & (f  \in \text{ span of CM-forms}) \\
     & X &  (\text{otherwise}).
     \end{aligned}
    \end{cases}
\end{align}
We say that $f\in $ span of CM-forms if it is a linear combination of Hecke eigenforms, each of whose underlying newform has complex multiplication by some (possibly different) imaginary quadratic fields. 
It is well known by a result of Serre that such an $f$ is \textit{lacunary} if $\kappa \ge 2$, cf. \cite[Theorem~17]{Serre81} and vice-versa. 
Recall that $f$ is called lacunary if $\displaystyle \# \{ n \le X\mid a_f(n) \neq 0\}=o(X)$ as $X \to \infty$. 
Thus, $f\in S_\kappa(\Gamma_1(N))$ is lacunary if and only if $\kappa=1$ (cf. \cite{Serre81}, \cite{DS74}) or $f$ is a CM-form for $\kappa\ge 2$. 
By abuse of notation we will use interchangeably the phrases ``$f$ is a CM-form" and that ``$f$ is in the span of CM-forms". 
Also it is clear that $f$ is CM if and only if all the newforms appearing in the newform-decomposition of $f$ (see \ref{se:int-wt-prilims}) are CM.

We now state our main result on elliptic cusp forms.  We make the usual convention that $4|N$ when $\kappa$ is a half-integer. We denote by  $\mathcal{U}_\kappa(N)$ the space of cusp forms spanned by the unary theta series (see \eqref{eq:ukn-defn}). This is relevant only when $\kappa=3/2$.

\begin{thm}\label{th:int-unrestricted}
    Let $N \ge 1$, $\kappa\in \frac{1}{2}\mathbb Z $ with $\kappa>0$, and $0 \neq f \in S_\kappa (\Gamma_1(N))$  be as in \eqref{expansion}. Then there exists $M=M_f\geq 2$ such that for any $\eps>0$ and sufficiently large $X,$
            \begin{align}
              \mathscr N_f(X;M) \gg \mathscr S_f (X;M) \gg_{f, \eps} 
\begin{cases}
    \begin{aligned}
        & \mathscr N_f (X), && \text{ if } f \in \mathcal{U}_\kappa(N),\\
        &X^{1-\eps}, &&  \text{ if } f \in S_\kappa(\Gamma_1(N))\setminus \mathcal{U}_\kappa(N), \\
        &\mathscr N_f (X), && \text{ if } \kappa \in \Z. 
    \end{aligned}
\end{cases} 
\end{align}
One can take $M=2$ when $\kappa$ is integral.
\end{thm}

The same conclusions extend, after a suitable rescaling, to cusp forms on arbitrary congruence subgroups. 
To the best of our knowledge, a result of this strength is presently not available in the literature for arbitrary cusp forms of general level and weights $\kappa \ge 1/2$.
In particular, if $N=1$, then \eqref{eq:serre-estimates}  shows that $\mathscr S_f(X) \asymp X$.
Even this result seems to be new, and can be thought of as the integral weight counterpart to \cite[Theorem 1.2]{mat15}. 

Theorem \ref{th:int-unrestricted} will follow as a consequence of Theorem \ref{thma} given below. 
Theorem \ref{thma}  will also play a crucial role in applications to Siegel modular forms.
We define the following hypothesis for the ease of presentation.

\noindent {\bf $\medcirc$ {\color{blue}Hypothesis~H:}}
\tagtext{Hypothesis~H}{tag:1}
There exists $n_0\ge 1$ with $(n_0,2N)=1$ such that $a_f(n_0)\neq 0$. 

Many of our arguments are easier to handle if $f$ satisfies \ref{tag:1}. 
Unlike the setting in \cite{jls}, where any non-zero $f$ in the Kohnen's plus space of level $4N$ with $N$ being odd, square-free automatically satisfied \ref{tag:1} by virtue of the results of Serre-Stark \cite{serre-stark}; such may not the case for general levels or integral weights. Further, we denote by  $ S_\kappa^* (\Gamma_1(N))$  the orthogonal complement of the space $\mathcal U_\kappa(N)$.

\begin{restatable}{thm}{HypH} \label{thma}
    Let $N \ge 1$, $\kappa\in \frac{1}{2}\mathbb Z $ with $\kappa \ge 1$, and $f \in S_\kappa^* (\Gamma_1(N))$  be as in \eqref{expansion} satisfying  \ref{tag:1}. 
    Then there exists $M\geq 2$ such that for any $\eps>0$ and sufficiently large $X,$
            \begin{align}
    \mathscr N_f^\flat(X;M) \gg \mathscr S_f^{\flat}(X;M) \gg_{\eps, f} 
    \begin{cases}
    \begin{aligned}
        &X^{1-\eps}, &&  \text{ if } \kappa \in \tfrac12 \Z,\\ 
        &\mathscr N_f^{\flat}(X), && \text{ if } \kappa \in \Z.
        \end{aligned}
\end{cases} 
\end{align} 
One can take $M=2$ when $\kappa$ is integral. 
Moreover when $\kappa \in \Z$, we have $\mathscr S_f^{\flat}(X) \asymp \mathscr N_f^\flat(X) \asymp \mathscr N_{f^\circ}(X)$, where $f^\circ$ is defined in \eqref{eq:f-circ}.
\end{restatable}

To place these results in context, we briefly recall what was known previously. 
The following table summarises the best known results in the literature. Below, by level $N$, we mean $\Gamma_0(N)$.

\begin{table}[ht]
\centering
\footnotesize
\setlength{\tabcolsep}{3pt}
\renewcommand{\arraystretch}{1.12}

\begin{tabularx}{\linewidth}{
@{}l>{\raggedright\arraybackslash}X
@{\hspace{1.2em}}
l>{\raggedright\arraybackslash}X@{}
}
\multicolumn{2}{c}{\textbf{Integral weight}} &
\multicolumn{2}{c}{\textbf{Half-integral weight}} \\

& $\mathscr{S}_f(X)$ &
& $\mathscr{S}^\flat_f(X)$ \\

\cmidrule(lr){1-2}\cmidrule(l){3-4}

\cite{mat15} &
$\gg_f X$; $f$ a Hecke eigenform of level $1$ &
\cite{jls} &
$\gg_{f,\varepsilon}X^{1-\varepsilon}$;
$f$ arbitrary in the Kohnen plus space,
level $4N$ ($N$ odd, square-free)
\\[2pt]

\cite{lau17} &
$\gg_f X^{1/2}$; $f$ and $N$ arbitrary 
 &
\cite{les21} &
$\gg_f X^{1-\eps}$; $f$ eigenform in the Kohnen's plus space level $4$\\

\cite{hul17} &
$\gg_f X^{2/5}$; $f$ and $N$ arbitrary 
 &
\cite{lau17} &
$\gg_f X^{1/2}$; $f$ and $N$ arbitrary
\end{tabularx}
\vspace{.2em}
\caption{Previous sign-change bounds for cusp forms.}
\label{tab:prior-sign-changes}
\end{table}
\vspace{-1.5em}
The main novelty of our results lies in their generality and quantitative strength: we cover the most general weights and levels, and obtain the best possible orders of magnitude of the quantities $\mathscr S_f^{\flat}(X;M)$ or $\mathscr S_f(X;M)$ in integral weight. 
Our results also improve upon the known lower bound  $\mathscr N_f^{\flat}(X) \gg X^{1-\eps}$ (\cite[Theorem 2]{ana19}, \cite[Theorem 4.5(b)]{bd22}) for an arbitrary cusp form (see Corollary \ref{cor:nonvanishing}).

Major advances in this topic have been possible by the use of sieve methods and by proving the desired results in suitable sets of  positive relative measure. We also employ an amalgam of both of these two ideas. In particular our methods are inspired by both \cite{mat15} and \cite{jls}, but significant further extensions and new challenges have to be overcome. 

For instance, we note that the method of \cite{jls} would more or less straightforwardly give the lower bound (at least on \ref{tag:1}) 
\[  \mathscr S_f^{\flat}(X;M) \gg_{\eps, f} X^{1-\eps}. \]
Therefore the entire efforts in the Sections~\ref{se:positive},~\ref{se:one-non-cm} and~\ref{se:all-cm} are to improve this lower bound to the exact correct order of magnitude of $\mathscr S_f^{\flat}(X;M)$. The distinction between CM and non-CM $f$ becomes especially prominent in this context. In addition, the small weights, viz. weights $\kappa=1/2,3/2$ and $\kappa=1$, require somewhat separate treatment, see Section~\ref{sec:32alt},~\ref{se:theta} and~\ref{se:all-cm}.
This occupies the major part of the paper. 

We can also consider the sign-change problem over primes. 
To this end, let us define \[\mathscr S_f^{\mathrm{pr}}(X) :=\#\{X< p_i \le 2X~:~p_i \text{ prime},\ p_i\nmid 2N,\ a_f(p_i)a_f(p_{i+1})<0\}.\] 
In Section \ref{prime}, we prove the following sign change result.

\begin{thm}\label{th:prime-intro}
Let $f\in S_\kappa(\Gamma_1(N))$ ($\kappa \in \Z, \kappa \ge 2$) be a cusp form satisfying \ref{tag:1} with the Fourier expansion \eqref{expansion}. Then there exists a constant $0<\vartheta_0<1$ depending only on $f$ such that for any $\eps>0$ and sufficiently large $X,$ we have $\displaystyle \mathscr S_f^{\mathrm{pr}}(X) \gg_{f,\eps} X^{1-\vartheta_0 - \eps}$.

Assuming Generalised Riemann hypothesis, we have $\mathscr S_f^{\mathrm{pr}}(X) \gg X(\log X)^{-(2+\eps)}$.
\end{thm}

\subsection{Brief idea about the proofs of the Theorems}
Let us sketch the main ideas behind the proofs. 
The proofs of Theorems \ref{thma} and Theorem \ref{th:prime-intro} are based on the standard principle commonly employed in the works dealing with sign changes: the strict inequality
\begin{align}\label{eq:sign-changes}
    \Big| \sumn_{n \in I} c_f(n) \Big| < \sumn_{n \in I} |c_f(n)|
\end{align} 
for suitable intervals $I$ enforces a sign change in $I$.
Improvements result by adding sieve weights to the above quantities.

The proof of the half-integral weight case of Theorem \ref{thma}, given in Section \ref{se:half-int}, follows broadly the argument in \cite{jls}, but there are notable differences. 
For instance, due to the general setting (non square-free levels), one does not immediately get access to central $L$-values via  Waldspurger's formula. 

In the case of the integral weights, we use `normalised' sieve weights $w(n)$ and $w'(n)$ \emph{away from the level} with $\displaystyle 0\leq w'(n)\leq w(n)$ (see \eqref{eq:wprime},~\eqref{eq:wupper}),  which allow us to exploit the `partial' multiplicativity inherited from linear combination of newforms (see Section \ref{se:positive} and \eqref{eq:bilinear}). We can then show that there exist many intervals $I$ with the property
\begin{align} \label{core-inequality}
    \Big|\sumn_{n\in I} \sgn(c_f(n))w(n) \Big| < \sumn_{n\in I} w'(n).
\end{align} 

The coefficients $c_f(n)$ are not multiplicative, which prevents a clean separation of `large' and `small' prime factors of $n$ as in \cite{mat15}. This feature is required for invoking the shifted convolution result from \cite{har03}.
Cancellation in the LHS of \eqref{core-inequality} is indeed detected by considering its second moment and reducing it to the bounds of suitable shifted-convolution problems. 
However, `large' factors in our sieve weights  hinder a direct application of known shifted-convolution estimates  uniformly in its parameters. See Lemma~\ref{lem:shift} and \eqref{eq:cdr}. We overcome this by truncating $\sgn(c_f(n))w(n)$ first (see \eqref{W-small}), and then applying \cite{har03} with a negligible error. This is inspired by \cite{jls}.
On the other hand, comparing the first and second moments of $w'(n)$ yields a lower bound for the RHS of \eqref{core-inequality}. 

Since lacunarity affects the natural size of the sieve-weight moments, one must separate the proof into two parts, depending on the CM-type of the newforms. 
The integral weight case of Theorem \ref{thma} is proved  by treating the non-lacunary and lacunary cases separately in Sections \ref{se:one-non-cm} and \ref{se:all-cm-11}, respectively. These are finally assembled in Section~\ref{sec:thma-end}.
This method heavily uses the `partial' multiplicativity, and thus cannot be used for the half-integral cases. 

In Section~\ref{se:theta}, we treat the weight $1/2$ case and the space of unary theta series of weight $3/2$ together. Here the argument is simple: since by Serre-Stark \cite{serre-stark}, everything is a linear combination of theta series, the sign changes can be detected in an interval of length the LCM of the respective conductors of the characters, by exploiting the orthogonality of characters.

The proof of Theorem \ref{th:prime-intro}, given in Section \ref{prime}, follows the same sign-detection principle (as for the weights $\kappa \ge 2$), applied to log-weighted sums over primes. 
The required cancellation and lower bounds follow from short-interval estimates for the coefficients of the relevant standard and Rankin-Selberg $L$-functions.

\subsection{Applications of Theorems \ref{th:int-unrestricted} and \ref{thma}}
We discuss some applications of our results on sign changes in the elliptic setting.
For instance it shows oscillation of the quantity $r_Q(n)-\mathrm{Eisen}_Q(n)$: the error term in representation numbers $r_Q(n)$ of a quadratic form $Q$ in $2k$ variables (cf. the remark in \cite[p.~1940]{mat15}). 

Moreover, let $Q_1,Q_2$ be positive-definite integral quadratic forms in $2k$
variables, $k\geq2$ is integral, belonging to the same genus, and assume that
$\theta_{Q_1}\ne\theta_{Q_2}$.  Write
\[
 D_{12}(X)=\#\{n\leq X:r_{Q_1}(n)\ne r_{Q_2}(n)\}.
\]
Since $\theta_{Q_1} - \theta_{Q_2}$ is a cusp form (\cite[Theorem~17]{Maass1954}), we immediately obtain that the quantities $r_{Q_1}(n)- r_{Q_2}(n)$, with $n \le X$, change signs at least $\gg D_{12}(X)$ times; and the same holds for the quantities $\{ n \le X \mid r_{Q_1}(n) \lessgtr r_{Q_2}(n)\}$. The same type of results hold e.g. for ternary quadratic forms, where we land in weight $3/2$.

\subsection{Results on holomorphic Siegel cusp forms}
For Siegel cusp forms of degree $n\ge 2$, Fourier coefficients indexed by fundamental discriminants form a natural arithmetically significant family which are of particular interest in automorphic representations. We refer the reader to \cite{sah13, sah-1-13, bd22, sha23} and \cite{jls} for more perspective on this topic. 

Hecke eigenvalues and Fourier coefficients are distinct objects for Siegel modular forms of degree $n \ge 2$. For a historical background and in particular results about the sign changes  related to eigenvalues, we refer the reader to \cite{pit08, das18, das20, pau23,jls} and the references therein. 
There are also interesting connections with mathematical physics. Fourier coefficients of certain meromorphic Siegel modular forms encode degeneracies of BPS states, with constrained sign patterns arising in the study of single-centred black holes; see \cite{sen2011black}.

In degree 2, \cite{jls} established $\gg_{F,\varepsilon}X^{1-\varepsilon} $ sign changes among fundamental Fourier coefficients with discriminant $\asymp X.$ 
Their method reduces the problem, via the Fourier-Jacobi expansion and the Eichler-Zagier correspondence, to sign changes of half-integral weight modular forms. 
The scalar-index structure in degree $2$ preserves the relevant congruence and (e.g. square-free) level conditions, which is crucial for obtaining the essentially optimal exponent $X^{1-\varepsilon}$. 

One of the purposes of the present paper is to extend such quantitative sign-change results to Siegel cusp forms of arbitrary degree. In higher degrees, the relevant (or the most convenient) Fourier-Jacobi coefficients generally have matrix indices, so the scalar-index reduction available in degree 2 no longer applies directly. Instead, using the theta decomposition and the structural results of B\"ocherer-Das \cite{bd22}, we isolate suitable primitive theta components and reduce the problem to Fourier coefficients of an elliptic cusp form. This route overcomes the lack of a general injectivity result for the corresponding Eichler-Zagier map, at the cost of changing the congruence subgroup and increasing the level. Moreover, it seems that the idea of using twisted Eichler-Zagier maps (cf. \cite{ana19, was26}) does not work here, as it is difficult to count sign changes from a situation in which complex valued characters and their linear combinations are involved.

However, we need very general results about sign changes in the elliptic case, since we lose control of both the congruence subgroup and its level in the process, as discussed above. This is one of the reasons that we have to work with full generality in the earlier sections.

We briefly recall the notation needed to state our main results. We define $\varrho(n) $ to be $ \frac n2+1$ or $\frac{n+3}{2}$ according as $n$ is odd or even, respectively.
Let $\Lambda_n^{+}$ denote the set of positive definite half-integral symmetric matrices of degree $n$. Further, for an $M \in \Lambda_n^+$, we define its absolute discriminant by \[d_M:= \begin{cases}
\det(2M), & n \text{ even},\\ \dfrac12\det(2M), & n \text{ odd}.
\end{cases}\]

 We write $\flat$ for the set of odd square-free positive integers and $\mathbb P$ for the set of primes and define 
\begin{align} \label{sft}
\mathscr S_{F}(X, \flat;M) := \#\left\{ d_j \in \flat \cap [X,MX]: \exists \ \mathcal T_j \in \Lambda_n^+,~ a_F(\mathcal T_j)a_F(\mathcal T_{j+1})<0,~d_j=d_{\mathcal T_j } \right\}.
\end{align} 
We will  drop $M$ and the subscript $F$ from the notation, as they will be evident from the context. 
Moreover, we define $\mathscr S(X,\mathbb P;M)$ analogously to $\mathscr S(X,\flat;M)$ in $\eqref{sft}$, with $\flat$ replaced by $\mathbb P$.

\begin{thm}\label{siegelpositiveprop}
Let $n\geq 2$  and $F\in S_k^n$ a non-zero Siegel cusp form  of weight $k\geq \varrho(n)$ with real Fourier coefficients. 
Then there exist infinitely many $T\in \Lambda_{n}^{+}$ with $d_T \in \flat$, and  constants $M\ge 2$, $X_0\ge 1$ (depending only on $F$ and $\eps$) such that for all $X\geq X_0$, 
\begin{align}
    \mathscr S(X,\flat)\gg \begin{cases}
        X \, (\log X)^{-1/2}   & \text{ if } n \text{ is odd}; \\
        X^{1-\eps} & \text{ if } n \text{ is even}.
    \end{cases}
\end{align}
When $n$ is odd, we moreover have $\displaystyle \mathscr S(X,\mathbb P)\gg X^{1-\delta-\eps}$ for some $0 < \delta <1$.
\end{thm}
The proof is given in Section~\ref{sec:smf}.
Here $\delta=\delta(F,T):=\vartheta_0(f),$ where $0<\vartheta_0(f)<1$ is the constant arising from Theorem \ref{unconditional:prime} applied to the elliptic cusp form $f$ associated with $F$ and $T$ as in \eqref{away} obtained from the theta expansion of a suitable Fourier-Jacobi coefficient of $F$. This is worked out in \cite{bd22}.
Under GRH, we can take $\delta=0$.

Note that, due to the distinct discriminant condition, the matrices $\mathcal T_j$ in the above theorem are $\GL_n(\Z)$ inequivalent. 

We also note that the above theorem sharpens the quantitative non-vanishing result  $\mathscr N(X,\flat)\gg_\eps X^{5/8-\eps}$ of \cite{bd22} for scalar-valued cuspidal Siegel modular forms of even degree.
However, our results do not cover non-cusp forms, unlike that in \cite{bd22} or \cite{bd18}. It will be interesting to investigate this further.

\section{Notation and preliminaries}\label{se:notation}
For a commutative ring $R$ with unity, we denote by $M(n;R)$ the set of $n\times n$ matrices with entries in $R$.
The general and special linear groups are denoted by $\GL(n,R)$ and $\SL(n,R)$, respectively. For a matrix $M$, its transpose is denoted by $M^t$. If $A$ and $B$ are matrices of compatible sizes, we use the notation $A[B]:=B^tAB.$ For a real symmetric matrix $A$, the notations $A\geq 0$ and $A>0$ means that $A$ is positive semi-definite and positive definite, respectively. Let
\begin{align} \label{lambdan}
    \Lambda_n := \left\{ S=(s_{ij})\in M(n;\mathbb Q): S=S^t,\quad s_{ii}\in\mathbb Z,\quad s_{ij}\in \tfrac12\mathbb Z,\quad S\geq 0 \right\},
\end{align}
and let $\Lambda_n^+$ denote the subset of positive definite matrices in $\Lambda_n$. 

For $n\geq 1$, $P^-(n)$ denotes its least prime factor, $\omega(n)$ the number of distinct prime divisors of $n$, and $\mu(n)$ the M\"obius function. We write $\ind_{n\in X}$ for the indicator of a set $X$, that is, 
\[\ind_{n\in X} = \begin{cases}
1, & n\in X,\\
0, & n\notin X.
\end{cases}\]
Unless stated otherwise, the letters $p$ and $q$ always denote prime numbers.

The symbol $\tau$ is used in different contexts throughout the paper; its meaning will always be clear from the context. 
Finally, we say that a property $\mathbf{P}$ holds for almost all $x\in [X,2X]$ if $\meas\{x\in [X,2X]:x\text{ satisfies } \mathbf{P}\}=X(1-o(1))$ as $X\rightarrow \infty$, where $\meas$ is the standard Lebesgue measure on $\R$. 
In Section \ref{sec:32alt}, we will also use the standard Haar measure $\mu:=\mu_{\mathrm{Haar}}$ on $\R_{>0}$ defined as $d\mu_{\mathrm{Haar}} = x^{-1}dx.$

 Unless otherwise indicated, all implied constants appearing in the various estimates in this paper may depend on $f$ and on any other fixed parameters, whose dependence will be indicated where relevant. 
 In particular, in any estimate involving an arbitrary $\eps>0$, the implied constant will  depend on $\eps$, which may vary from line to line. 
 Thus we will simply write $A \ll B$ for $A \ll_{f,\eps} B$.

\subsection{Elliptic modular forms}\label{se:ell-prilims}
Throughout the paper we make the standing assumption that:
\begin{center}
    $\kappa \ge  \frac 12$ is allowed to be either integral or half-integral, and

$N \ge 1$ is arbitrary if $\kappa$ is integral; and is a multiple of $4$ if $\kappa$ is half-integral.
\end{center}
We briefly recall the notation and basic facts concerning elliptic modular forms that will be used throughout the paper. For integral and half-integral weight modular forms, we refer the reader to the standard references \cite{dia05} and \cite{shi73}, respectively. 

Let us denote by $S_\kappa(\Gamma)$ the inner product space of cusp forms for a congruence subgroup $\Gamma$ of $\SLtwoZ$. We will be mainly concerned with cusp forms for the congruence subgroup $\Gamma_1(N)$. Recall that $\Gamma_0(N)$ denotes the set of matrices $\left(\begin{smallmatrix}
    a&b\\c&d
\end{smallmatrix}\right)$ in $\SLtwoZ$ such that $N\mid c$, and $\Gamma_1(N)$ is the subset of $\Gamma_0(N)$ of matrices with $a\equiv d\equiv 1 \pmod{N}.$ We have the following decomposition of the vector space $S_\kappa(\Gamma_1(N))$ as a direct sum \[S_\kappa(\Gamma_1(N)) = \oplus_{\chi} S_\kappa(\Gamma_0(N),\chi),\]  where the sum runs over all Dirichlet characters $\chi$ $\!\!\!\!\mod N.$ For brevity, we write $S_\kappa(N,\chi)=S_\kappa(\Gamma_0(N),\chi) $.

Each $f\in S_\kappa(\Gamma_1(N))$ has a Fourier expansion of the form \eqref{expansion}.
We denote by $c_f(n)$ the `normalised' $n$-th Fourier coefficients of $f$ given by \[a_f(n) = c_f(n) n^{\frac{\kappa-1}{2}}.\] 

\subsubsection{Half-integral weight}\label{se:half-int-wt-prelims}
We first give some details on the theta subspaces that occur in weights $\frac12$ and $\frac32$. Put $\nu=0$ or $1$ according as $\kappa=\frac12$ or $\frac32$. 
For a primitive Dirichlet character $\psi$ of conductor $r_\psi$ with $\psi(-1)=(-1)^\nu$ and an integer
$t\ge1$, define 
\begin{align}\label{eq:defn-theta}
     \theta_{\psi,t}^{(\nu)}(\tau) =\sum_{m\ge1}m^\nu\psi(m)q^{tm^2}.
\end{align}
 The transformation formula for
theta series gives
\[
 \theta_{\psi,t}^{(0)}\in M_{1/2}(4r_\psi^2t,\psi\chi_t),  \qquad  \theta_{\psi,t}^{(1)}\in S_{3/2}(4r_\psi^2t, \psi\chi_{-4t}),
\] 
where $\chi_{t}$ is the usual Kronecker symbol.
For $\kappa=1/2$, the theta series $\theta_{\psi,t}^{(0)}$ is cuspidal precisely when $\psi$ is even but not \emph{totally even} (see \cite{serre-stark}). We define
\begin{align}\label{eq:ukn-defn}
    \mathcal U_{\kappa}(N) := \mathrm{span} \left\{ \theta_{\psi,t}^{(\nu)} ~:~    \psi(-1)=(-1)^\nu,~ \psi\text{ primitive, }  4r_\psi^2t\mid N     \right\}\cap S_{\kappa}(\Gamma_1(N)), \quad S_{\kappa}^*(N) = \mathcal U_{\kappa}(N)^\perp.
\end{align}
By \cite[Theorem~A]{serre-stark}, all cusp forms of weight $\frac12$ are linear combinations of such unary theta series, i.e.,
\[\mathcal U_{1/2}(N)=S_{1/2}(\Gamma_1(N)).\]
For weight $\kappa\ge 5/2$, there is no unary theta series, and we put $S_\kappa^*(\Gamma_1(N)) = S_\kappa (\Gamma_1(N)).$ 

\subsubsection{Integral weight}\label{se:int-wt-prilims}
For $\kappa\in \Z$, by the standard newform theory, every  $f\in S_{\kappa}(\Gamma_1(N))$ admits a decomposition of the form 
\begin{equation}\label{eq:newform-decomp-simplified}
  f(\tau)=\sumn_{1\le j\le r}\sumn_{d\mid (N/N_j)}\alpha_{j,d}\,f_j(d\tau)
\end{equation}
for some $r$ (see \cite[\S 5, Theorem 5.8.3]{dia05}), where $f_j$ are primitive forms of level $N_j\mid N$ which may have complex Fourier coefficients. 
Let us further define a cusp form
\begin{align}\label{eq:f-circ}
    f^\circ(\tau) =  \sumn_{1\le j\le s} \alpha_{j,1} ~f_j(\tau).
\end{align} 
We will always write $\alpha_{j}=\alpha_{j,1}$ for simplicity.
Then for every $(n,N)=1$, we have 
\begin{align}\label{eq:newform-decomp-simplified-coeffs}
    c_f(n) = c_{f^\circ}(n) = \sumn_{1\le j\le s} \alpha_{j}~ \lambda_j(n),
\end{align}
where $\lambda_j(n)$ are the normalised Fourier coefficients of the eigenforms $f_j.$

For every $p|N$ and a cusp form $h\in S_\kappa(\Gamma_1(N))$, define the $p$-sifting operator $\mathcal{D}_p$ by 
\begin{align}\label{eq:dp-operator}
    \mathcal D_p h := h-V_pU_p h,
\end{align} where $U_t$ and $V_t$ are defined as \[U_t \Big(\sumn_{n\ge 1} A(n) q^n\Big) =  \sumn_{n\ge 1} A(nt) q^n \quad \text{ and } \quad V_t \Big( \sumn_{n\ge 1} A(n) q^n \Big)=  \sumn_{n\ge 1} A(n) q^{nt}.\]  

\subsection{Counting sign changes}
Here we prove an elementary lemma for counting the sign changes, which will be invoked in the later sections.
\begin{lem}\label{le:sign-change-count}
    Let us call an interval $I_{x,y}=(x,x+y]$ to be  good if it contains a sign change of a real sequence $(a_n)_{n \ge 1}$, i.e.,
    \begin{align}
        |\sumn_{\substack{x< n\le x+y}}~ a_n| < \sumn_{\substack{x< n\le x+y}}~ |a_n|.
    \end{align}
   Put $\mathcal{G}_X:= \{x\in [X,MX]: I_{x,y} \text{ is good}\}.$ 
   Then 
    \begin{align}
        \mathscr{S}(X,M):= \#\text{sign~changes~of $a_n$ in~$[X,MX]$~} \gg \meas(\mathcal{G}_X)/y
    \end{align}
    whenever $\meas(\mathcal{G}_X)>y$.
    \end{lem}
\begin{proof}
    Let $\mathcal{S}_X$ denote the set of pairs $(m,n)$ of consecutive indices of non-zero  $(a_n)_{n\geq 1}$ with $a_ma_n<0$,  $X\leq m<n\leq MX$. Then, $\#\mathcal{S}_X = \mathscr{S}(X,M).$
    
    For every $x\in\mathcal{G}_X$ with $x\leq MX-y$, there exists a pair $(m,n) \in \mathcal{S}_X$ such that  $x< m< n \leq x+y$.
    In particular, $x \in [n-y,m)$. Thus \[ \mathcal{G}_X \subseteq (MX-y,MX] \cup \bigcup_{(m,n) \in \mathcal{S}_X} [n-y,m). \] 
    Taking measures on both sides, we obtain $ \meas(\mathcal{G}_X)\leq y(1+\#\mathcal{S}_X).$
    Since $\meas(\mathcal{G}_X) > y$, the claim now follows.
\end{proof}

\section{ \texorpdfstring{Proof of $\mathscr S_f^{\flat}(X;M) \gg X^{1-\eps}$ in Theorem \ref{thma} for  $\kappa\in \frac{1}{2}\Z$}{κ ∈ 12 Z}}\label{se:half-int}
Let $N \ge 1$, $\kappa\in \frac 12 \Z$ such that $\kappa \ge\frac{3}{2}$. In this section, we will prove that if $f\in S_\kappa^* (\Gamma_1(N))$ satisfies \ref{tag:1}, then
\begin{align} \label{k-half}
\mathscr{S}^{\flat}_f (X;M) \gg_{f,\eps} X^{1-\eps},
\end{align}
for some $M \ge 2$ depending only on $f$.
Our  aim is  to treat both integral and half-integral cases simultaneously whenever feasible. 

Note that the results of this section are known when $f$ belongs to Kohnen's plus space and has odd square-free level, see \cite[Sections 2, 3]{jls}. In this case the \ref{tag:1} holds automatically from the main results of Serre-Stark \cite[Lemma~7]{serre-stark}. Our aim is to show how to modify some of those arguments to accommodate arbitrary forms. This section is along the lines of \cite{jls}.

Before proceeding with the proof of the theorem, we establish the following auxiliary lemma, which will play a key role in the argument. 
\begin{lem}\label{propRS}
     Let $N\ge 1$, $\kappa \in \frac12 \Z$ with $\kappa \ge \frac{3}{2}$ and $f\in S_\kappa^* (\Gamma_1(N))$ be a cusp form with the Fourier expansion \eqref{expansion} satisfying $\ref{tag:1}.$ Then for some $M\geq 2$ and sufficiently large $X$, we have 
     \begin{equation} \label{eqRS}
  \sum_{X \leq n \leq MX , \, (n,2N)=1 } |c_f(n)|^2 \mu^2(n) \asymp_{f, M} X,
   \end{equation}
and for any $\eps>0$,
   \begin{equation} \label{fouthmoment}
         \sum_{X \leq n \leq 2X , \, (n,2N)=1}|c_f(n)|^4 \mu^2(n) \ll_{f, \eps} X^{1+\eps}.
   \end{equation}
   \end{lem}

\begin{proof}
    For integral weights, \eqref{eqRS} follows from the classical Rankin-Selberg method after sieving out squares; see \cite[Theorem 4.3, eq.~(4.25)]{bd22} for proofs. The fourth moment in \eqref{fouthmoment} follows directly from Deligne's bound.
     
   Now we consider the half-integral weight case.
   Under \ref{tag:1}, the bound \eqref{eqRS} for $\kappa > 3/2$ can be deduced from the arguments put together in \cite[Proposition~2.3, eq.~(14)]{jls}. Indeed, the upper bound in \cite[eq.~(16)]{jls} clearly holds for all $f \in S_\kappa^*(\Gamma_1(N))$ and all $N$ by the Rankin-Selberg asymptotic bound. When $\kappa \ge 5/2$, the lower bound in \eqref{eqRS} follows from \cite[\S 4.1.2, eq.~(4.31)]{bd22}, which adapts the argument of \cite[Proposition~3.7]{sah13} from $\Gamma_0(N)$ to $\Gamma_1(N)$.
   For $f \in S_{3/2}^*(\Gamma_1(N))$, in view of \cite[Proposition~2.4, Remark~2.5--Theorem~2.6]{li17}, the lower bound follows by the same argument as in \cite[\S 4.1.2]{bd22}, with the necessary modifications.

   The upper bound in \eqref{fouthmoment} requires more care. Without loss we can assume that $f\in S_\kappa^*(N,\chi)$, $\kappa\ge 3/2$, and $N$, $\chi$ are both arbitrary. In \cite{jls}, an explicit version of Waldspurger's formula was utilized. The required estimate was reduced to a Lindel\"of on average bound for the second moment of twisted central $L$-values $L(1/2, \mathbf f \otimes \chi_d)$, where $\mathbf f$ is the Shimura lift of $f$. When $N$ is not square-free,  we argue as follows.
   
   The adelisation $\widetilde\varphi=t(f)$ belongs to a cuspidal automorphic representation $\widetilde\pi$ of $\widetilde{\mathrm{SL}}_2(\mathbb A_{\mathbb Q})$. Without loss, we can assume that  $\widetilde\varphi=t(f)$ belongs to an irreducible cuspidal automorphic representation $\pi$. Let $g$ be the unique newform of level $\leq N^2$ and trivial nebentypus associated to $\pi$. Then, from  \cite[Appendix 2, \S~7.4]{blo07}, for a square-free positive integer $D,$ we have $$|c_f(D)|^2 \ll_{f,N} L(\tfrac{1}{2}, g\otimes \chi_D),$$ where $\chi_D$ is the primitive quadratic Dirichlet character associated with $\mathbb Q(\sqrt{D}).$ Consequently, we have 
   \begin{align}
       \sum_{X\le n\le 2X,\, (n,2N)=1} |c_f(n)|^4 \mu^2(n) \ll_{f} \sum_{|D| \le 2X,\, D ~\mathrm{fundamental}} |L(\tfrac{1}{2}, g\otimes \chi_D)|^2 \ll_{f,\eps} X^{1+\eps}
   \end{align} 
   which is a standard application of Heath-Brown's quadratic large sieve \cite{bro95}, see e.g. \cite[Corollary~2.5]{sou10}. Note that for us the level is fixed.
\end{proof}
\begin{rmk}
    In the integral weight case, clearly one can take $M=2$ in the above Lemma.
\end{rmk}

\subsection{Proof of the bound in \protect\eqref{k-half} -- Theorem~\ref{thma} for half-integral weights \texorpdfstring{$>$}{>} 3/2}\label{se:signintegral}
In this section we prove the lower bound in \eqref{k-half}. 
We first note two auxiliary results, viz. Propositions \ref{squarerootcancellation} and~\ref{absolutelowerbd} given below, which follow by suitably modifying the corresponding arguments given in \cite{jls}.
We next observe that once these two Propositions are in place; the proof of \eqref{k-half} is a formal verbatim deduction, and thus will be omitted. 
We choose to provide brief arguments for the other supplementary assertions for the sake of continuity. 
Our claim is that the results in this subsection hold in much generality than that is available in \cite{jls}, in particular for all $\kappa\in \frac{1}{2}\Z,$ $\kappa\ge 2.$ 
The first Proposition says that the mass of the square-free Fourier coefficients is $X^{-\eps}$ on average in most short intervals.
\begin{prop} \label{absolutelowerbd}
 For a positive integer $N$ and $\kappa\in \frac 12 \Z$ with $\kappa \ge \frac{3}{2}$, let $f \in S_\kappa^*(\Gamma_1(N))$ satisfy \ref{tag:1}. Then, for any $\eps >0$ and $2 \le y \le X/2$ there are $\gg_{f,\eps} X^{1-\frac32 \eps}$ integers $X \le x \le MX$ such that
 \begin{align}\label{eq:absolutelowerbd}
       \sum_{x \leq n \leq x + y , \, (n,2N)=1} |c_f(n)| \mu^2(n)  >  \frac{y}{X^{\eps}}.
 \end{align}
 \end{prop}
\begin{proof}
    When $\kappa\in\Z$, from \eqref{eqRS} and Deligne's bound, it immediately follows that 
    \begin{align}
        \sumn_{X \leq n \leq 2X , \, (n,2N)=1} |c_f(n)| \mu^2(n)  \gg_f X^{1-\eps/2}.
    \end{align}

    When $\kappa\in\frac{1}{2}\Z\setminus \Z$ (as in \cite{jls}), we apply H\"older inequality to $\sum_{\substack{X\le n\le MX\\(n,2N)=1}} |c_f(n)|^2 \mu^2(n)$ in terms of first and fourth moments. Then, using the estimates in Lemma \ref{propRS}, we obtain 
    \begin{align}
        \sumn_{X \leq n \leq MX , \, (n,2N)=1} |c_f(n)| \mu^2(n)  \gg_f X^{1-\eps/2}.
    \end{align}
Now proceeding exactly as in \cite{jls}, the proposition follows.
\end{proof}
 
 Besides showing that the absolute values of the Fourier coefficients are typically large in short intervals, one also needs sufficient cancellation in the corresponding signed sums. The next proposition provides such a cancellation result on average over short intervals.
\begin{prop}\label{squarerootcancellation}
     Let $N\ge 1$, $\kappa \in \frac12 \Z$ with $\kappa>\frac{3}{2}$ and $f \in S_{\kappa}(\Gamma_1(N))$.
     Then for $1 \leq y \leq X^{\frac{1}{209}}$ we have \[  \sum_{X \le x \le 2X} \bigg|\sum_{\substack{x \le n \le x+y  \\ (n,2N)=1}} c_f(n) \mu^2(n) \bigg| \ll_f X \sqrt{y}.\]
\end{prop}

\subsection{Proof of Proposition \ref{squarerootcancellation}}\label{se:sqroot}

Here, we briefly indicate the  changes necessary (as compared to the proof given in \cite[Proposition~3.3]{jls}) for proving the proposition. Several modifications are required to accommodate the more general situation considered here. We will need the following Lemma which hold for $\kappa=3/2$ also. 

\begin{lem}\label{lem:sumbd}
Let $N\ge 1$, $\kappa \in \frac12 \Z$ with $\kappa\ge \frac{3}{2}$ and $f \in S_{\kappa}^*(N,\chi)$. Then \[\sum_{n \le X, \ \ d^2\mid n ,\ \  (n,2N)=1  } |c_f(n)| \ll_{f,\eps} \frac{X^{1+\eps}}{d^2}.\]
\end{lem}
\begin{proof}
    We omit the arguments for $\kappa$ integral which follows directly from Deligne's bound.

    We give a proof when $\kappa$ is not integral and $f$ is not necessarily in the Kohnen's plus space. Since it is sufficient to prove the statement for eigenforms, we assume that $f$ is an eigenfunction for all Hecke operators $T_{p^2}$, where $p$ is a prime not dividing $N$, with the normalised eigenvalues $\lambda_f(p).$ i.e., \[T_{p^2} (f) = p^{k-\frac 12} \lambda_f(p) f.\]
    Recall that for a non-zero half-integral weight eigenform $f$, there exists a square-free integer $t$ such that its Shimura lift $g:=\mathrm{Sh}_t(f)$  is non-zero \cite{shi73}. 
    Moreover, $g\in S_{2\kappa-1}(N/2,\chi^2)$ is a normalised cuspidal Hecke eigenform 
    whose eigenvalues coincide with those of $f$; see \cite[Proposition 3.2.4]{pur12} for more details. In particular, for every prime $p\nmid N$, we have $\lambda_g(p)=\lambda_f(p).$  Thus from Deligne's bound, we immediately have $|\lambda_f(p)|\leq 2.$ Using this and the Hecke relations \cite[Corollary 1.8]{shi73} for half-integral weight forms \[c_f( p^{2j+2} n) = \lambda_f(p) c_f( p^{2j} n) - \chi(p^2)c_f ( p^{2j-2} n)\] for primes $p$ not dividing $N$, by an induction argument on $j$ we obtain \[|c_f(p^{2j}n)|\leq 3^{j} |c_f(n)|\] for all $j\geq 0$. 
    Note that by Cauchy-Schwarz and \eqref{eqRS}, $\sum\limits_{m\leq Y,~(m,2N)=1}|c_f(m)|\mu^2(m) \ll_f Y.$ Now, writing $n=t^2m$ with $m$ square-free, we have \begin{align*}
        \sumn_{\substack{n \le X,~  (n,2N)=1 \\ d^2 |n }} |c_f(n)| \ll_\eps X^{\eps} \sumn_{\substack{t \le \sqrt{X},~ d|t}} ~ \sumn_{ m \le \frac{X}{t^2} } |c_f(m)|  \ll_\eps \frac{X^{1+\eps}}{d^2},
    \end{align*} as desired.
\end{proof}
The following proposition holds for both integral and half-integral weight cusp forms almost verbatim from \cite[Proposition 6.1]{les21}. 
\begin{prop}\label{shifted}
    Let $N\ge 1$, $\kappa \in \frac12 \Z$ with $\kappa>\frac{3}{2}$ and $f\in S_\kappa(N,\chi)$  with the Fourier expansion \eqref{expansion}. Then for every $\eps>0$, the following estimate holds uniformly for all integers $1 \leq \Delta \leq X^{1/52}$, $0<|h|<X^{1/2}$ and $(v,\Delta)=1$.
\begin{equation}\label{shiftedconv}
\frac{1}{X^{\kappa-1}} \sumn_{n\geq 1} c_f(n) n^{\frac{\kappa-1}{2}} e\!\left(\frac{nv}{\Delta} \right) e^{-2\pi n/X} c_f(n+h) (n+h)^{\frac{\kappa-1}{2}} e^{-2\pi(n+h)/X} \ll X^{1-\frac{1}{104}+\eps}.
\end{equation}
\end{prop}
\begin{proof}
    The proof follows the application of Jutila's version of the circle method as presented in \cite[Proposition~6.1]{les21}, with the choice of moduli used in \cite[Proposition~3.5]{jls}.
    We observe that it holds for general levels and any $\kappa \ge \frac52$.
    Namely, we take $q=N\Delta p$, where $p\equiv1\pmod4$ is prime, and apply \cite[Lemma~6.1]{les21} with its parameter $4\Delta$ replaced by $N\Delta$.     
\end{proof}
We now give a very brief sketch of the proof of Proposition \ref{squarerootcancellation}.
\begin{proof}[\textbf{Proof of Proposition \ref{squarerootcancellation}}]
We closely follow the method of \cite[Proposition 3.3]{jls}. We write $\mu^2(n) =\sum_{d^2|n} \mu(d)$ and define
\[\mu^2_{\leq Y}(n) := \sumn_{\substack{d^2 |n,~d \leq Y }} \mu(d)  \qquad  \mu^2_{> Y}(n) := \sumn_{\substack{d^2 |n,~ d > Y }} \mu(d).\]
Then by Lemma \ref{lem:sumbd}, we have the contribution from the `large' divisors is at most \begin{equation} \label{eq:largedbd}
\sum_{X \leq x \leq 2X} \bigg|\sum_{\substack{x \leq n \leq x + y \\ (n,2N)=1}} c_f(n) \mu^2_{>Y }(n) \bigg| \leq y \sum_{\substack{n \leq 4 X \\ (n,2N)=1}} |c_f(n) \mu^2_{> Y}(n)| \ll y \sum_{\substack{d > Y }} \sum_{\substack{n \leq 4 X \\ d^2 | n}} |c_f(n) |  \ll \frac{y X^{1+\eps}}{Y}.
\end{equation}

For the `small' divisors, we denote by $A_f(n) = c_f(n) \mu^2_{\leq Y }(n) \mathbf 1_{(n,2N)=1},$ and consider the weight function $W(t) = t^{\frac{\kappa-1}{2}} e^{-2\pi t}.$ Then a similar argument as in \cite{jls} gives that \begin{equation}\label{shift!}
    \sum_{X \le x \le 2 X} \Big | \sum_{x \le n \le x+y} A_f(n) \Big|\ll \sqrt{X}  \Big( \sum_{x \ge 1} \Big| \sum_{x \le n \le x + y } A_f(n)  \Big |^2 W\left( \frac{x}{X}\right)^2\Big)^{1/2}.
\end{equation}
The infinite sum inside the parenthesis equals (cf. \cite[Equation (25)]{jls})
\begin{equation}\label{shifthh}
    \sumn_{|h| \leq y} (y + 1 - |h|) \sumn_{\substack{n \ge 1 }} A_f(n) A_f(n + h) W \Big (\frac{n}{X} \Big) W \Big ( \frac{n + h}{X}  \Big ) + O_{f, \eps}\left(y^3  X^{\eps} \right).
\end{equation}
The diagonal contribution of $
\ll_{f, \eps}~ yX$ follows from \eqref{eqRS}. We follow the arguments of \cite{jls} to obtain that the contribution from $h\neq 0$ as 
\[\ll y   \sum_{\substack{0 < |h| \leq y \\ d_1, d_2 \leq Y  \\ d_3,d_4 | 2N}}  \bigg| \frac{1}{r} \sum_{0 \le v < r} e\left(\frac{-av}{r} \right) \sum_{n \ge 1} c_f(n) c_f(n+h) e \left ( \frac{n v}{r} \right ) W \Big ( \frac{n}{X} \Big ) W \Big ( \frac{n + h}{X} \Big )\bigg|\]
with $r \leq 16N^2 Y^4.$ 
Applying Proposition \ref{shifted} to the inner sum (for $(v,r)>1$, we reduce the fraction to lowest terms), we see that \eqref{shifthh} is  
\begin{align}\label{eq:est-half-final}
    \ll yX+ y^2 Y^2 X^{1-\frac{1}{104}+\eps} + y^3 X^{\eps},
\end{align} 
provided that $r\leq 16N^2 Y^4 \leq X^{\frac{1}{52}}$ and $|h|< \sqrt{X}.$ 
Consequently, \eqref{shifthh}  is $\ll yX$ whenever, in addition, $  y Y^2 X^\varepsilon \le X^{1/104}$ and $y^2 X^\varepsilon\le X.$ With $Y=\sqrt y X^\varepsilon$, these conditions become 
\[16 N^2 y^2 X^{4\varepsilon}\le X^{1/52},  \qquad y^2 X^{3\varepsilon}\le X^{1/104}, \qquad  y^2 X^\varepsilon \le X.\]
They hold for $1\le y\le X^{1/209}$ and sufficiently large $X$, upon choosing $\varepsilon$ sufficiently small. 
It follows from \eqref{shift!} and \eqref{eq:est-half-final} that the small-divisor contribution is $\ll X\sqrt y$. 
Finally, from \eqref{eq:largedbd}, for $Y=\sqrt y\,X^\varepsilon$, large-divisor contribution is also $\ll \sqrt y\,X$, which completes the proof.
\end{proof}
As we have discussed in Section~\ref{se:signintegral}, this finishes the proof of \eqref{k-half} -- Theorem~\ref{thma} for half-integral weights.
\begin{rmk}
    Probably the results of the above Section holds also for $\kappa=3/2$ but we have not checked this fully. The only difference could be in the shifted convolution estimate, which is stated for weights $\kappa \ge 5/2$ in \cite{jls}.
\end{rmk}

\subsection{An alternative proof of Theorem \ref{thma} for \texorpdfstring{$\kappa\ge 3/2$}{κ ≥ 32 }} \label{sec:32alt}
In this section we give a proof of Theorem \ref{thma} by a mean-square argument of $L$-series attached to $f$, avoiding the shifted convolution argument of \cite{jls}. 
To this end, unlike in the previous sections, we will work with `multiplicatively' short intervals of the form $(x,xe^\delta]$, $0<\delta<1$, since these are particularly well suited to the application of Plancherel (see \cite{gal70, sou17}).

As the integral weight case is similar and simpler, we only demonstrate the non-integral weight case. 
Let $F\in S_{\kappa}(L,\chi)$ and $\mathcal{L}(s,F) = \sumn_{n\ge 1} c_F(n) n^{-s}$. 
The completed function $\Lambda_F(s)= (\sqrt L/2\pi)^{s+(\kappa-1)/{2}} \Gamma(s+\tfrac{\kappa-1}{2}) \mathcal L_F(s)$ is entire and satisfies $\Lambda_F(s) = \omega_F \Lambda_{F^*}(1-s)$, where $F^*(\tau) = (\sqrt L \tau)^{-\kappa} F(-1/(L\tau))$ is the Fricke transform of $F$ and $|\omega_F|=1$. 
Note that $B_{F^*}=B_F$, where \[B_F:=\sup_{\tau\in\mathbb H} (\Im\tau)^{\kappa/2} |F(\tau)|.\]
Further, we have the following approximate functional equation (see \cite[Theorem 5.3]{iwa04}):
\begin{align}\label{eq:afe}
    \mathcal{L}(s,F) = \sum_{n\ge1} \frac{c_F(n)}{n^s} V_s(n) + \varepsilon(s,F) \sum_{n\ge1} \frac{c_{F^*}(n)}{n^{1-s}}V_{1-s}(n), \quad  |\varepsilon(1/2+it,F)|=1,
\end{align}
where $V_s(u)$  is a smooth function defined by
\begin{align}\label{eq:afe-weight}
    V_s(u) = \frac{1}{2\pi i} \int_{(3)} e^{w^2}\frac{\Gamma(s+\tfrac{\kappa-1}{2}+w)}{\Gamma(s+\tfrac{\kappa-1}{2})}\Big(\tfrac{\sqrt L}{2\pi u}\Big)^{w} \frac{dw}{w}, \quad 0\le \Re(s)\le 1.
\end{align}
Moreover, for any $U\ge 2$ with $U/2\le |t| \le 3U$,  $x>0$, from \eqref{eq:afe-weight} and Stirling's estimates applied to the ratio of the Gamma functions, it follows that (\cite[Proposition 5.4]{iwa04}) 
\begin{align}\label{eq:afe-weight-bound}
    |V_{1/2+it}(x)| + x|V_{1/2+it}'(x)| \ll (1+x/Q_0)^{-1},
\end{align}
where $Q_0 \asymp_\kappa \sqrt{L}U$, the implied constant does not depend on $F$. 
We will also use the following  mean-value theorem for  Dirichlet polynomials.
\begin{lem}[{\cite[Corollary 9.2]{iwa04}}]\label{le:mont-mean-value}
For any complex coefficients $b_n$, we have 
    \begin{align}
        \int_0^{T} \left|\sumn_{1\le n\le N} b_n n^{-it} \right|^2 dt \ll (T+N)\sumn_{n\le N} |b_n|^2.
    \end{align}
\end{lem}
Using this, we will deduce the following uniform mean-square estimate of $L(\tfrac{1}{2} + it, F)$.
\begin{lem}\label{le:mean-value-l}
    Let $F\in S_\kappa(L,\chi)$ with $4\mid L$. Then for $U\ge 2$, we have 
    \begin{align}\label{eq:mean-value-l}
        \int_U^{2U}|\mathcal L(\tfrac{1}{2} + it, F)|^2 dt \ll B_F^2 U(\sqrt L+\log^2 2U).
    \end{align}
\end{lem}
\begin{proof}
    The proof is standard but a bit subtle in the present case since we need explicit control on the implied constants. Later $L$ will be replaced with $L=L_0r^2$ with $L_0$ fixed. We therefore give enough details for the convenience of the reader. Put $\mathcal{P}_F(s):=\sum_{n\ge1} \frac{c_F(n)}{n^s} V_s(n) $. 
    Let 
    \begin{align}
        \mathcal P_M(t) = \sum_{M\le n<2M} c_F(n) n^{-1/2-it} V_{1/2+it}(n), \quad \text{and} \quad \mathcal R_x(t) =  \sum_{M\le n< x} c_F(n) n^{-1/2-it}.
    \end{align}
    Let us denote by $\|\cdot\|_2$ the $L^2$-norm taken over $[U,2U]$.
    From Lemma \ref{le:mont-mean-value} and Parseval, we have 
\begin{align}\label{eq:rx2}
     \|\mathcal{R}_x\|_2 \ll (U+M)^{1/2} \Big( \sumn_{M\le n\le 2M} \frac{|c_F(n)|^2}{n}\Big)^{1/2} \ll B_F (U+M)^{1/2}, \quad M\le x\le 2M,
\end{align}
Now, for  fixed $t$, Abel summation gives
    \begin{align}\label{eq:pmt-rxt-1}
        \mathcal P_M(t) = V_{1/2+it}(2M)~ \mathcal R_{2M}(t) - \int_M^{2M} V_{1/2+it}'(x) ~\mathcal R_x(t)dx.
    \end{align}
Taking $L^2$-norm on both sides above, and then using Minkowski's integral inequality 
\[\left( \int_U^{2U} \left| \int_Y f(t,y) dy \right|^2 dt \right)^{1/2} \le \int_Y \left( \int_U^{2U} |f(t,y)|^2 dt \right)^{1/2} dy\] 
for $L^2$-spaces in the last term and \eqref{eq:afe-weight-bound}, we get $\| \mathcal{P}_M \|_2  \ll B_F {\sqrt{(U+M)}}/{(1+M/Q_0)}.$
Since $\mathcal{P}_F(1/2+it) = \sumn_{M\text{ dyadic}} \mathcal P_M(t)$, splitting the sum into $M\le U$, $U<M\le Q_0$, and $M>\max\{U,Q_0\}$ gives 
\begin{align}\label{eq:afe-l2-bound-1st-sum}
    \|\mathcal{P}_F(1/2+it)\|_2 \ll B_F( \sqrt{U} \log (2U) + \sqrt{Q_0}).
\end{align}
The same estimate holds for the dual sum in \eqref{eq:afe} after replacing $t$ with $-t$. 
Since $Q_0 \asymp_\kappa U\sqrt{L}$, combining \eqref{eq:afe-l2-bound-1st-sum} and the corresponding estimate for the dual sum with \eqref{eq:afe}, we finally have \eqref{eq:mean-value-l}.
\end{proof}

Let us fix $0< \delta \le 1$ and $K>1$ throughout the section.
\begin{prop}\label{pr:mean-square-logarithmic}
    Let $F\in S_\kappa(L,\chi)$ and $S_F(x,\delta):=\sumn_{x< n\le xe^\delta}c_F(n).$
    Then for $Y\ge 1$ and $\eps>0$,  
    \begin{align}\label{eq:mean-square-logarithmic}
        \int_Y^{2Y} |S_F(x,\delta)|^2 \frac{dx}{x} \ll_{\eps} B_F^2\delta^{1-\eps} Y \sqrt{L}.
    \end{align}
\end{prop}
\begin{proof}
    Choose a smooth function $W\in C_c^\infty(0,\infty)$ such that $W(u)=1$ for $u\in [1,2e]$. Further, for $v\in \R$, define 
    \begin{align}
        \Phi(v) := \sumn_{n\ge 1} c_F(n) W(n/Y) \ind_{v<\log n \le v+\delta}, \quad \text{and} \quad D_Y(t)=\sum_{n\ge1}c_F(n)W(n/Y)n^{-it}.
    \end{align}
    Then for $Y\le x\le 2Y$, we have  $S_F(x,\delta) = \Phi(\log x)$ and $\int_Y^{2Y} |S_F(x,\delta)|^2 \frac{dx}{x} \le \int_\R |\Phi(v)|^2 dv.$
    Note that Mellin inversion gives
    \begin{align}
        D_Y(t) = \frac{1}{2\pi i} \int_\sigma \Tilde{W}(s) Y^s \mathcal{L}(s+it,F) ds = \frac{\sqrt{Y}}{2\pi} \int_\R \Tilde{W}(1/2+iv)Y^{iv} \mathcal L(1/2+i(t+v),F) dv,
    \end{align}
    where $\Tilde{W}(s)$ is the Mellin transform of $W$ and $\sigma > 1$. Now applying Minkowski's inequality we have 
    \begin{align}\label{eq:dyt-l2}
        \int_{-V}^V |D_Y(t)|^2 dt \ll B_F^2 YV(\sqrt{L} + \log^2(2V)) \ll V^{1+\eps} Y \sqrt{L}.
    \end{align}
    On the other hand, from Plancherel's formula it follows that 
    \begin{align}\label{eq:plancherel-32}
        \int_\R |\Phi(v)|^2 dv &=\frac{1}{2\pi} \int_\R |D_Y(t)|^2 \left| \frac{e^{it\delta}-1}{it}\right|^2 dt.
    \end{align}
    Write the right-hand side as $\int_{|t|\le \delta^{-1}}  + \int_{|t|> \delta^{-1}} $.
    We use the pointwise bound $|(e^{it\delta}-1)/(it)| \le  \delta$ and \eqref{eq:dyt-l2} in the first integral to obtain
    \begin{align}\label{eq:le-delta}
        \int_{|t|\le \delta^{-1}}  |D_Y(t)|^2 \left| \frac{e^{it\delta}-1}{it}\right|^2 \ll B_F^2 Y\delta^{1-\eps} \sqrt{L}.
    \end{align}
    Similarly, in the complementary region, we use $|(e^{it\delta}-1)/(it)| \le 2/|t|$ and get
    \begin{align}\label{eq:ge-delta}
         \sumn_{j\ge 0} \int_{2^j/\delta < |t| \le 2^{j+1}/\delta} |D_Y(t)|^2 \left| \frac{e^{it\delta}-1}{it}\right|^2 \ll \delta^2 \sumn_{j\ge 0} 2^{-2j} \int_{ |t| \le 2^{j+1}/\delta} |D_Y(t)|^2 dt \ll B_F^2 Y\delta^{1-\eps} \sqrt{L}.
    \end{align}
    Now, combining \eqref{eq:plancherel-32} and the above estimates, we have the assertion.
\end{proof}

 For the remaining part of the section, let us fix  
 \begin{align}\label{eq:delta-defn}
     y=X^{17 \eps} \quad \text{and} \quad \delta=y X^{-1},
 \end{align} and let $X\ge X_0(f,K,\eps)$ be large enough that $2\le y\le X/2$.  Put for any $K>1$,
    \begin{align}\label{eq:bx-ex-defn}
        \mathcal{B}(x) = \sumn_{\substack{x< n\le xe^\delta\\ (n,2N)=1}} c_f(n) \mu^2(n), \quad \text{and} \quad \mathcal E_X := \Big\{x\in [X,KX]~:~ |  \mathcal{B}(x) | \ge \tfrac{1}{3}yX^{-\eps} \Big\}. 
    \end{align}
\begin{prop}\label{pr:signed-32}
    For a fixed $K>1$ and $f\in S_\kappa^*(\Gamma_1(N))$ satisfying \ref{tag:1}, we have $\mu(\mathcal E_X) \ll_{f,K,\eps} X^{-2\eps}$.
\end{prop}
\begin{proof}
    We start by defining \[h(\tau) := \sumn_{\substack{n\ge 1\\ (n,2N)=1}} a_f(n) q^n,\] which is a non-zero cusp form (since $f$ satisfies \ref{tag:1}) and has level $N^3$ (\cite[\S 4.2, Lemma 4]{serre-stark}).
    For some fixed $R$ (to be determined later),  write $\mathcal{B} = \mathcal{B}_{\le R} + \mathcal{B}_{>R},$    where
    \begin{align}\label{eq:brx-description}
         \mathcal{B}_{\le R}(x) = \sumn_{r\le R} \mu(r) \sumn_{x<r^2 m \le xe^\delta} c_h(r^2m),
    \end{align}
    and     $\mathcal B_{>R}$ is defined likewise with $r>R$. Here we have used the identity $\mu^2(n) = \sumn_{r^2|n} \mu(r)$.

     For each $r\le R$, define $h_r := h\mid_\kappa U_{r^2}$. Then $h_r\in S_\kappa^*(\Gamma_1(N^3 r^2))$ (cf. Lemma \ref{le:duke-bound}) and $c_h(r^2 m) = r^{-(\kappa-1)} c_{h_r}(m).$ 
     Therefore, from \eqref{eq:brx-description} we have 
     \begin{align}\label{eq:brx-hr}
         \mathcal{B}_{\le R}(x) = \sumn_{r\le R} \mu(r) r^{-(\kappa-1)} S_{h_r}(x/r^2,\delta).
     \end{align}
     By Cauchy-Schwarz we have 
     \begin{align}
         |\mathcal{B}_{\le R}(x)|^2 \le  \sumn_{r\le R } \mu^2(r) \sumn_{r\le R }   |r^{-2(\kappa-1)} S_{h_r}(x/r^2,\delta)|^2.
     \end{align}
     We have the {\it uniform in} $r$ bound $ \frac{B_{h_r}}{r^\kappa} \le B_h \ll_f 1$ which follows by bounding the summands in $h\mid_\kappa U_{r^2}$ absolutely.
     Now, Proposition \ref{pr:mean-square-logarithmic} applied to $h_r$ with $Y\asymp X/r^2$, $L=N^3 r^2$ and $B_{h_r}^2\ll_f r^{2\kappa}$ gives (for some $C_f>0$)
     \begin{align}
         \int_X^{KX} |\mathcal{B}_{\le R}(x)|^2 \frac{dx}{x} &\le R\sumn_{r\le R} r^{-2(\kappa-1)} \int_{X/r^2}^{KX/r^2} | S_{h_r}(u,\delta) |^2 \frac{du}{u} \\ 
         &\ll_{K,\eps}  R \sumn_{r\le R} B_{h_r}^2 (X/r^2) \delta^{1-\eps}  \frac{\sqrt{N^3 r^2} }{r^{2(\kappa-1)}}   
         \le C_f^2 X \delta^{1-\eps} R^{3}  \ll_{f,K,\eps} R^3 y X^\eps.
     \end{align}
     
     Also, by the triangle inequality and Lemma \ref{lem:sumbd},
     \begin{align}
         \int_X^{KX} |\mathcal{B}_{> R}(x)| \frac{dx}{x} 
         \ll \delta  \sumn_{r>R} \sumn_{\substack{n\le KXe^\delta\\r^2|n,~ (n,2N)=1}} |c_f(n)| 
         \ll_{f,K,\eps} \delta  \sumn_{r>R} \frac{X^{1+\eps}}{r^2} \ll \frac{yX^\eps}{R}.
     \end{align}
     Now, for any $x\in \mathcal{E}_X$, we must have at least one of  $|\mathcal{B}_{\le R}(x)|$ and $|\mathcal{B}_{> R}(x)|$ is larger than $\frac{1}{6}yX^{-\eps}$. 
     Let $\mathcal{C}_X:=\{x\in [X,KX]~:~|\mathcal{B}_{\le R}(x)| \ge \frac{1}{6}yX^{-\eps}\}$ and  $\mathcal{D}_X$  the corresponding set for $\mathcal{B}_{> R}(x)$. 
     Then from Chebyshev's and Markov's inequality respectively, we have 
     \begin{align}
         \mu(\mathcal{C}_X) \ll \frac{R^3 X^{3\eps}}{y}, \quad \text{and} \quad \mu(\mathcal{D}_X) \ll \frac{X^{2\eps}}R.
     \end{align}
     Taking $R=X^{4\eps}$ gives $\mu(\mathcal E_X)\le\mu(\mathcal C_X)+\mu(\mathcal D_X)\ll X^{-2\eps}$.
\end{proof}

\subsubsection{Proof of Theorem~\ref{thma} for half-integral weights}\label{se:thma-mean}
By Proposition \ref{absolutelowerbd}, we already have a set $\mathcal{V}_X\subset [X,KX]$ ($K\ge M$) of integers $m$ with  $|\mathcal{V}_X| \gg_{f,\eps} X^{1-3\eps/2}$ such that $\sumn_{\substack{m \le n \le m+h\\(n,2N)=1}} |c_f(n)|\mu^2(n) > hX^{-\eps}$ holds.
Let us define 
\begin{align}
    \mathcal{W}_X:= \Big\{ x\in [X,KX]~:~  \sumn_{\substack{x \le n \le xe^\delta\\(n,2N)=1}} |c_f(n)|\mu^2(n) > \frac{1}{3}yX^{-\eps} \Big\}.
\end{align}
Note that for each $m\in \mathcal{V}_X$ and   $h=\lfloor y\rfloor-1$, we have $(m-1,m)\cap [X,KX] \subset \mathcal{W}_X$. Therefore, \[\mu(\mathcal{W}_X) \ge \sumn_{m\in \mathcal{V}_X} \mu((m-1,m))\gg_{f, \eps} \frac{X^{1-3\eps/2}}{KX} \gg_{f,K,\eps} X^{-3\eps/2}.\]
Recall the definition of $\mathcal{E}_X$ from \eqref{eq:bx-ex-defn}. 
For every $x\in \mathcal{W}_X \setminus \mathcal{E}_X$, we have 
\begin{align}
    \sumn_{\substack{x \le n \le xe^\delta\\(n,2N)=1}} |c_f(n)|\mu^2(n) > \left| \sumn_{\substack{x< n\le xe^\delta\\ (n,2N)=1}} c_f(n) \mu^2(n) \right|,
\end{align}
and from Proposition \ref{pr:signed-32}, $\mu(\mathcal{W}_X \setminus \mathcal{E}_X)\ge \mu(\mathcal{W}_X) - \mu(\mathcal{E}_X)\gg X^{-3\eps/2}$. 
Since $\meas(A) \ge X\cdot \mu(A)$ for any set $A\subset [X,\infty)$, and a good `sign-detecting' interval $(x,xe^\delta]$ has length $x(e^\delta-1)\asymp_K X\delta =y$ for any $x\in [X,KX]$, from Lemma \ref{le:sign-change-count}, we have \[\mathscr{S}_f^\flat(X;K) \gg \frac{X^{1-3\eps/2}}{y} \gg_{f,K,\eps} X^{1-\frac{37\eps}{2}},\] and we are done. \qed
\section{Theta series}\label{se:theta}
In this section we will consider the sign changes among the Fourier coefficients of cusp forms of weights $\kappa=\frac 12, \frac 32$ defined by unary theta series  as in \eqref{eq:defn-theta}.
Note that the theta series can be written as 
\begin{align}
     \theta_{\psi,t}^{(\nu)}(\tau) = \sumn_{ s\equiv 0 (\!\bmod{d}),~ t=t_0d^2}  (s/d)^\nu\psi(s/d) ~q^{t_0s^2}.
\end{align}
Therefore, for any  $f= \sumn_{\psi, t} \alpha_{\psi,t} \theta_{\psi,t}\in \mathcal{U}_\kappa(N)$ and $n=t_0 s^2$, $t_0$ square-free, the $n$-th Fourier coefficient of $f$ is given by 
\begin{align}\label{eq:afn-theta-combination}
    a_f(t_0 s^2) = s^\nu P_{t_0}(s),\qquad P_{t_0}(s) = \sum_{d\mid s} \sum_\psi \alpha_{\psi,t_0d^2} d^{-\nu} \psi(s/d).
\end{align}
Let $\mathscr T$ be the finite set of square-free $t_0$ with $P_{t_0} \not\equiv 0$. For $t_0 \in \mathscr T$, we claim 
\begin{enumerate}
    \item $a_f(t_0s^2) \neq 0 \iff P_{t_0}(s)\neq 0$, and the two have the same sign,
    \item $ P_{t_0}$ is periodic with period $M_{t_0} = \lcm\{dr_\psi : \alpha_{\psi, t_0 d^2} \neq 0 \}$, and 
    \item $\sumn_{1\le s\le  M_{t_0}} P_{t_0} (s) = 0$.
\end{enumerate}
The claim (i) is obvious from the definition. 
For (ii), note that for any $d\ge 1$ and non-principal $\psi$ with conductor $r_\psi$, the function $F_{d,\psi} (s)$ defined by 
\begin{align}
    F_{d,\psi} (s) = \begin{cases}
        \psi(s/d),  &\text{if } d\mid s,\\
        0, &\text{otherwise},
    \end{cases}
\end{align}
is periodic with period $dr_\psi$.\\ For $t_0\in \mathscr{T}$, let $A_{t_0} := \{(d,\psi) :  \alpha_{\psi,t_0d^2} \neq 0\}$. Then $P_{t_0}(s) = \sumn_{(d,\psi)\in A_{t_0}} \alpha_{\psi,t_0d^2} d^{-\nu} F_{d,\psi}(s)$, and consequently, $P_{t_0}(s) = P_{t_0}(s+M_{t_0})$.

Also, \[\sumn_{1\le s\le M_{t_0}} P_{t_0}(s) = \sumn_{(d,\psi)\in A_{t_0}} \alpha_{\psi,t_0 d^2} \sumn_{1\le s\le M_{t_0}} F_{d,\psi}(s).\]
The innermost sum is $ \frac{M_{t_0}}{d r_\psi} \sumn_{\substack{1\le s\le d r_\psi\\ d\mid s}} \psi(s/d) = \frac{M_{t_0}}{d r_\psi} \sumn_{1\le m\le r_\psi} \psi(m)$. But $\psi$ being non-principal, the last sum vanishes. This proves (iii).

In particular, $P_{t}(s)$  changes sign among $1\le s\le M_{t}$ for each $t\in \mathscr{T}$. 
Note that for each such $t$, one has $\mathscr S_f(X) \ge \#\{\text{sign changes of } P_t(s), \sqrt{X/t} \le s\le  \sqrt{2X/t} \}.$ 
Therefore, fixing $t\in\mathscr T$ and counting only $n=ts^2\in [X,2X]$,
\begin{align}\label{eq:nf-sf-theta-1}
\mathscr{N}_f(X) \ge \mathscr S_f(X) &\ge \left\lfloor \frac{1}{M_t} \left(\sqrt{\frac {2X}{t}} - \sqrt{\frac Xt}\right) \right\rfloor \gg_f \sqrt X.
\end{align}
for any $t\in \mathscr{T}$ and sufficiently large $X$. 

Moreover, since the square classes are pairwise disjoint,  from the periodicity of $P_t$ it follows that 
\begin{align}\label{eq:nf-theta-1}
    \mathscr{N}_f(X) = \sumn_{t\in \mathscr{T}} \#\{\sqrt{X/t}\le s\le \sqrt{2X/t}~:~ P_t(s)\neq 0\} =  \sqrt{X} \cdot \sumn_{t\in \mathscr{T}} \tfrac{K_t(\sqrt{2}-1)}{M_t\sqrt{t}}  + O_f(1)  \asymp \sqrt{X},
\end{align}
where $K_t = \#\{1\le s\le M_t : a_f(ts^2)\neq 0\}. $
From \eqref{eq:nf-sf-theta-1} and \eqref{eq:nf-theta-1}, the claim of Theorem \ref{th:int-unrestricted} now follows.

\section{Proof of Theorem~\ref{thma} when \texorpdfstring{$\kappa\in\Z$}{κ ∈ Z}: Setting up the sieve } 
\label{se:positive}

Let $f\in S_\kappa(\GammaOne{N})$ be a cusp form with Fourier expansion \eqref{expansion} satisfying \ref{tag:1}. We introduce non-negative sieve weights $w'(n) \leq w(n)$ and compare
\begin{equation}\label{eq:wnwn'sums}
\Big | \sum_{x \leq n \leq x + h} \sgn(c_f(n)) w(n) \Big | \quad \text{ and } \quad \sum_{x \leq n \leq x + h} w'(n).
\end{equation}
A strict inequality between these quantities forces a sign change in $[x,x+h]$. 
This approach leads to a proof of Theorem~\ref{thma} in the case of integral weights.  
In this section we will set-up the sieve machinery required to attain this objective. 

Due to the natural appearance of the CM forms and their sparsity of the  non-vanishing coefficients, our proof of Theorem~\ref{thma} for integral weights (thus assuming \ref{tag:1}), will be divided into two parts. 
\begin{itemize}
    \item Section \ref{se:one-non-cm} will treat the case that $f^\circ$ (as in \eqref{eq:f-circ}) is non-lacunary. 
    \item Section \ref{se:all-cm} will treat the complementary case. 
    \item Section~\ref{sec:thma-end} will assemble these together to finish the proof.
\end{itemize}
We begin by defining the weights $w(n)$ and $w'(n)$, and the notation used throughout the argument. 

Let $X\ge 1$ be fixed, and choose $y=X^{\delta},$ where $0<\delta < \frac{1}{2}$ is taken sufficiently small.
Following \cite[Section 3]{mat15}, we define $y_m = y^{\frac{1}{2}(1-\gamma^2) \gamma^{m-1}}$ with a parameter $\gamma \in (0, 1)$ and 
\begin{equation*}
 \mathcal{D}^+=\{d = p_1 \cdots p_\ell : p_1> \cdots > p_\ell, p_m \le y_m\text{ for every odd }m\} \subseteq[1,y], \end{equation*}
where the primes are distinct, and the empty product $1$ is included. Put \[\rho^+(b) = \sumn_{\substack{e\mid b,~e\in \mathcal{D}^+}} \mu(e).\]

The following standard upper-bound Brun-sieve properties were used in \cite[Sections 3 and 5]{mat15}.
\begin{align}
&0 \le \ind_{P^-(b)>y} \le \rho^+(b) \le \tau(b),\label{eq:sieve-positive}\\
&\rho^+(b) \le \ind_{P^-(b)>y} + \sumn_{\ell \ge 0} 4^{-\ell} \ind_{P^-(b)> y_{2\ell+1}} 4^{\omega(b)}.\label{eq:sieve-remainder}
\end{align}

With the decomposition \eqref{eq:newform-decomp-simplified-coeffs}, for $Q=2N$ and $(a,Q)=1$, let
\begin{equation}\label{eq:V} 
V(a)^2=\sumn_{1\le i \le s} |\alpha_i\lambda_i(a)|^2. \end{equation} Now define
\begin{equation}\label{eq:wprime} 
    w'(n)=\mu^2(n) \ind_{(n,Q)=1} \sum_{\substack{ab=n\\a\le y,\ V(a)>0}} \ind_{P^-(b)>y} \frac{|c_f(n)|}{V(a)}, 
\end{equation}
and
\begin{equation}\label{eq:wupper} 
    w(n)=\mu^2(n) \ind_{(n,Q)=1} \sum_{\substack{ab=n\\a\le y,\ V(a)>0}}\rho^+(b) \frac{|c_f(n)|}{V(a)}. 
\end{equation}

By~\eqref{eq:sieve-positive}, we have $ 0\leq w'(n) \leq w(n), $ with both weights vanishing whenever $c_f(n)=0$.

Since $c_f(n)$ is not multiplicative in general, we cannot simply discard the small factors $a$ as in \cite{mat15}. Instead we normalise by $V$, writing 
\begin{equation}\label{eq:beta} 
\beta_i(a) := \frac{ \alpha_i \lambda_i(a)}{V(a)}. \end{equation}
Then $\sumn_{1\le j\le s} |\beta_i(a)|^2 = 1$ (in particular $|\beta_i(a)| \le 1$), and multiplicativity of each $f_i$ gives
\begin{equation}\label{eq:bilinear}
\frac{c_f(ab)}{V(a)} = \sumn_i \beta_i(a) \lambda_i(b) \qquad \bigl((ab,Q)=1,\ (a,b)=1,\ V(a)>0 \bigr). 
\end{equation}
This lets us retain small factors without letting their size interfere with the sieve. The sieve weights can then be controlled using only the eigenvalues at the large factors, which allows us to apply the Brun-sieve majorant and the mixed shifted-convolution estimates of \cite{har03} uniformly in $a$.

Our goal is to compare the sums in \eqref{eq:wnwn'sums} with $h$ large but \textit{fixed} positive constant. 
When $f^\circ$ is lacunary, the short intervals requires an additional factor of $(\log X)^{-\theta_\circ}$ in their length; see Section \ref{se:all-cm}. 
We conclude the section by the following lemma which will be crucially used in Sections \ref{se:one-non-cm} and \ref{se:all-cm}.
\begin{lem}\label{le:first-moment-w'-preli}
 Let $\kappa\in \Z$ with $\kappa\ge 1$, and $f\in S_\kappa(\GammaOne{N})$   with Fourier expansion \eqref{expansion} satisfying \ref{tag:1}. Then for fixed $0<A<B$ and sufficiently large $X$, we have
    \begin{align}\label{eq:lower-bound-w'-preli}
        \sumn_{AX <n \le BX} w'(n) &\gg_f X(\log X)^{-1} H_V(y); \quad \text{where}\quad H_V(y):= \sumn_{\substack{a\leq y,\ (a,Q)=1\\ V(a)>0}}~ \frac{\mu^2 (a)}{a}.
    \end{align}
\end{lem}
\begin{proof}
    The prime Rankin-Selberg theorem  (see, for example, \cite[Section 5, (5.5)]{bd22})  for the fixed primitive forms $f_i$ and $f_j$ gives
    \[    \sumn_{\substack{At<q\le Bt,~ q\nmid Q}} \lambda_i(q) \conj{\lambda_j(q)} \log q = (B-A)\delta_{i,j} t + o_f(t). \]
    For $\beta=(\beta_i)_{1\le i \le s}$ with $\sum_{1\le i\le s}|\beta_i|^2=1$, we have $\sumn_{\substack{At<q\le Bt,~q\nmid Q}}|\sum_i\beta_i\lambda_i(q)|^2\log q = (B-A)t+o_f(t).$
    Now, Deligne's bound trivially gives $|\sum_{1\le i\le s} \beta_i \lambda_i(q)| \le 2 s$. Thus, using $\log q\le\log(2t)$, we obtain \begin{equation}\label{eq:prime-average}
        \sumn_{\substack{At <q \le Bt,~ q\nmid Q}} \Big| \sumn_i \beta_i \lambda_i(q) \Big| \ge \frac{B-A}{4s} \frac{t}{\log t} \gg_f\frac{t}{\log t},
    \end{equation} 
    for sufficiently large $t$. 

For $aq\in [X,2X]$ and $a\le y = X^\delta$, we have $q> X/y > y\ge a$ (since $\delta < 1/2$). Consequently $(a,q)=1$, $P^-(q)>y$, and if $a$ is square-free then $aq$ is square-free. 
Retaining just these factorisations in \eqref{eq:wprime} and using \eqref{eq:bilinear}, we have
\begin{align}\label{eq:first-lower}
    \sum_{AX < n \le BX} w'(n) &\ge \sum_{\substack{a\le y,~ (a,Q)=1,~ V(a)>0}}  \mu^2(a) ~\sum_{\substack{AX/a < q \le BX/a,\ q\nmid Q}} \left| \sumn_i \beta_i \lambda_i(q) \right|\\ &\gg_f X \sum_{\substack{a\le y,~(a,Q)=1,~ V(a)>0}} \frac{\mu^2(a)}{a \log(X/a)}\gg_f \frac{X}{\log X} H_V(y),
\end{align} 
completing the proof.
\end{proof}

\section{\texorpdfstring{Proof of Theorem~\ref{thma} when \texorpdfstring{$\kappa\in\Z$}{κ ∈ Z}: $f^\circ$}{f ◦ } is non-lacunary}\label{se:one-non-cm}
In this section, we will show that  
\begin{align} \label{asymp-nonCM}
    \mathscr{S}^{\flat}_f(X) \asymp \mathscr{N}_f^\flat(X) \asymp X
\end{align}
under the assumption that $f^\circ$ is non-lacunary, i.e., some $f_j$ in \eqref{eq:newform-decomp-simplified-coeffs} is non-CM. 

\subsection{Lower bound for \texorpdfstring{$w'$}{w′}}\label{se:lower-bound-w'}
We first establish a lower bound for $H_V(y)$ defined in \eqref{eq:lower-bound-w'-preli}   for later use.
\begin{lem}\label{le:noncm}
    Let $f^\circ$ be non-lacunary. Then we have $ H_V(y) \gg_f \log y.$
\end{lem}
\begin{proof}
Since $V(a) = \Big(\sumn_i |\alpha_i\lambda_i(a)|^2\Big)^{1/2} \ge |\alpha_i\lambda_i(a)|$ for any fixed $f_i$, one has $V(a)>0$ whenever $\lambda_i(a)\ne 0$. So it suffices to prove, for one fixed non-CM newform $f_i$ with $\alpha_i\ne0$ (without loss let's assume $j=1$), 
\begin{align}\label{eq:hvy-single}
    \sumn_{\substack{a\le y,~ (a,Q)=1\\ \lambda_{f_1}(a)\ne 0}} \frac{\mu^2(a)}{a}\ \gg\ \log y.
\end{align}
Now, let  $\mathcal M=\{m \ge 1: p\mid m, p\nmid Q \Rightarrow p^2\mid m\}$. 
For $n\le y$ with $\lambda_{f_1}(n) \neq 0,$ we can write $n=am$ with $a:=\prod_{\substack{p||n\\p\nmid Q}} p$. 
Then  $m\in\mathcal M$ and $\lambda_{f_1}(a) \ne0$, and we have
\begin{align}\label{eq:hvy-compo}
    \sum_{n \leq y,~\lambda_{f_1}(n) \neq 0}\frac{1}{n}~\le~ \sum_{\substack{a\le y,~ (a,Q)=1\\ \lambda_{f_1}(a) \ne 0}}\frac{ \mu^2(a)}{a}\ \sum_{m\in\mathcal M}\frac1m = C_\mathcal{M} \sum_{\substack{a\le y,~ (a,Q)=1\\ \lambda_{f_1}(a) \ne 0}}\frac{ \mu^2(a)}{a}.
\end{align}
By unique factorisation $C_\mathcal{M} := \sum_{m\in\mathcal M}\frac 1 m = \prod_{p\mid Q} \Big(1-\frac1p\Big)^{-1}  \prod_{p\nmid Q}\Big(1+\frac{1}{p(p-1)}\Big) <\infty.$ 
Moreover, from \cite{Serre81}, $\mathscr{A}(t):=\#\{1\le n\le t:  \lambda_{f_1}(n)\ne0\}\gg t$.
Hence, by partial summation,  the leftmost sum in \eqref{eq:hvy-compo} is $\gg_{f_1} \log y$, and \eqref{eq:hvy-single} follows.
\end{proof}

\begin{prop}\label{pr:moments}
    If $f^\circ$ is non-lacunary, then there exists a constant $c_0>0$ (depending only on $f$) such that 
    \begin{equation}\label{eq:main-first} 
        \sumn_{X< n\le2X}~w'(n) \ge c_0 X.
    \end{equation}
        Moreover, the following uniform pointwise bound holds.
        \begin{equation}\label{eq:main-first-pointwise}
            \sup\nolimits_{X< n\le 3X}~w'(n)\ll_f 1.
        \end{equation}
\end{prop}

\begin{proof}
The first assertion follows immediately from Lemma \ref{le:first-moment-w'-preli}  by putting $y=X^\delta$ in Lemma \ref{le:noncm}.

To prove \eqref{eq:main-first-pointwise}, recall that for $(a,Q)=1$, we have $V(a)^2 = \sum_{1 \le j \le s} |\alpha_i \lambda_i(a)|^2. $ 
Note that from multiplicativity and Cauchy-Schwarz,
\begin{align}
    |c_f(ab)| &={} \Big| \sumn_i \alpha_i \lambda_i(a) \lambda_i(b) \Big| \leq \Big(\sumn_i \Big|\alpha_i \lambda_i(a) \Big|^2\Big)^{1/2}  \Big(\sumn_i \Big|\lambda_i(b) \Big|^2\Big)^{1/2}, \nonumber
\end{align}
and therefore \[\frac{|c_f(ab)|}{V(a)} \leq \Big(\sumn_i \Big|\lambda_i(b) \Big|^2\Big)^{1/2}.\] 
Also, for a square-free $b$, we have $|\lambda_i(b)|\leq 2^{\omega(b)}.$ 
Consequently, $w'(n) \leq (s\ 2^{2\omega(b)})^{1/2} =2^{\omega(b)} \sqrt{s}. $ 
Since $b\leq 3X$ and every prime factor of $b$ is strictly greater than $y=X^\delta,$ we have $y^{\omega(b)}\leq b \leq 3X$, i.e.,
\[\omega(b) < \frac{3\log X}{\log y}= O_\delta(1).\] 
Consequently, from \eqref{eq:wprime}, we conclude that $w'(n) \ll_{f,\delta} 1$ and $\eqref{eq:main-first-pointwise}$ follows.
\end{proof}

\subsection{Upper bound for \texorpdfstring{$w$}{w}}
We use the following mean-value bound for multiplicative functions. 
\begin{lem}[{\cite[Theorem III.3.5]{Tenenbaum95}}]
\label{le:MVmultfunct}
Let $g \colon \mathbb{N} \to [0, \infty)$ be a multiplicative function. Let $A, B$ be constants such that 
for all $y \geq 1$,
\[\sum_{p \leq y} g(p) \log p \leq Ay \quad \text{and} \quad \sum_{p\in \mathbb{P}} \sum_{\nu \geq 2} \frac{g(p^\nu)}{p^\nu} \log p^\nu \leq B.\]

Then, for $x \geq 1$,
\[\frac{1}{x} \sum_{n \leq x} g(n) \ll (A+B+1) \prod_{p \leq x} \left(1-\frac{1}{p} + \sum_{\nu \geq 1} \frac{g(p^\nu)}{p^\nu}\right),\]
where the implied constant is absolute.
\end{lem}
The mean-value estimate now gives the required second-moment bound. 
Note that the following estimate holds regardless of whether $f^\circ$ is lacunary or non-lacunary. In the lacunary case, we will establish a stronger upper bound in the next section.

\begin{prop}\label{pr:upperbound2}
    For sufficiently large $X,$ we have \[\sumn_{n\le 3X}w(n)^2\ll_f X.\]
\end{prop}

\begin{proof}
Since $\rho^+ \ge 0$, we can drop the restrictions on $a$ involving $V(a)$ and $Q$ in \eqref{eq:wupper} to obtain
\[w(n) \le \sqrt{s} \mu^2(n) \sumn_{\substack{ab=n,~a\le y}} ~\rho^+(b)2^{\omega(b)}. \] 
Inserting \eqref{eq:sieve-remainder} in the above equation, we obtain that 
\begin{equation}\label{eq:w-majorant1} 
    w(n) \le \sqrt{s} \mu^2(n) \sumn_{\substack{ab=n,~a\le y}} \Big( \ind_{P^-(b)>y} + \sumn_{\ell\ge0} 4^{-\ell} \ind_{P^-(b)>y_{2\ell+1}} 2^{2 \omega(b)}\Big) 2^{\omega(b)}. 
\end{equation}

In its first term, $P^-(b)>y$ and $a\le y$ force $a$ to contain precisely the prime factors of $n$ at most $y$. Thus there is at most one such factorisation and therefore the first term is just $O_\delta(1).$

For the remainder put $z_\ell=y_{2\ell+1}$. After multiplication by $2^{\omega(b)}$, its $\ell$-th contribution is at most 
\[ 4^{-\ell}\sum_{\substack{ ab=n,\,a\le y\\P^-(b)>z_\ell}} 8^{\omega(b)} \leq  ~4^{-\ell} \sum_{\substack{ab=n\\P^-(b)>z_\ell}} 8^{\omega(b)} \le ~4^{-\ell}  \prod_{\substack{q\mid n\\q>z_\ell}} (1+8) = ~4^{-\ell} ~9^{\omega_{>z_\ell}(n)},\]
where $\omega_{>z_\ell}(n) = \#\{q\mid n~:~q>z_\ell\}$. 
Therefore
\begin{equation}\label{eq:w-nine} 
    w(n) \ll_{\delta,s} \mu^2(n) \Big(1 + \sumn_{\ell\ge0} 4^{-\ell} 9^{\omega_{>z_\ell}(n)} \Big). 
\end{equation}
Since $\sum_{\ell\ge0}4^{-\ell}<\infty$, Cauchy-Schwarz gives
\begin{equation}\label{eq:w-eightyone} 
    w(n)^2 \ll_{\delta,s} \mu^2(n) \Big( 1 + \sumn_{\ell\ge0} 4^{-\ell} 81^{\omega_{>z_\ell}(n)} \Big). 
\end{equation}

For $z>1$, define the non-negative multiplicative function \[g_z(n):=\mu^2(n)81^{\omega_{>z}(n)}.\] 
Then, for a prime $q$, we have
\[g_z(q^\nu)=
\begin{cases}
1, & \text{if } \nu=1 \text{ and } q\le z,\\
81, & \text{if } \nu=1 \text{ and } q>z,\\
0, & \text{if } \nu\ge 2.
\end{cases}\]
Thus, the mean-value upper bound in Lemma \ref{le:MVmultfunct} applies uniformly in $z$ with $A$ sufficiently large constant and $B=0$, and consequently, 
\begin{align}\label{eq:gz}
\sumn_{n\le 3X} g_z(n) &\ll X \prod\nolimits_{q\le 3X} \left(1 - \frac1q + \frac{g_z(q)}q \right) \ll X \begin{cases}
    \left( \frac{\log(6X)}{\log z} \right)^{80}, & \text{ for } z \geq 2,\\
    \log^{80}(6X) & \text{ for } 1 < z < 2.
\end{cases}
\end{align}
For $z\ge2$, the last step follows from Mertens' estimate 
\begin{align}\label{eq:merten}
\sumn_{z < q \le 3X} q^{-1} = \log \log (3X) - \log \log z + O(1)
\end{align} 
and $\log(1+t) \le t$. 
For $1<z<2$, we have used the product over all primes up to $3X$.

Recall that $y_m = y^{\frac{1}{2}(1-\gamma^2) \gamma^{m-1}}$, $y=X^\delta$, with a fixed  parameter $\gamma \in (0, 1)$ and  $z_\ell=y_{2\ell+1}$.  If $z_\ell\ge2$, then \eqref{eq:gz} and $\log(6X)\ll\log X$ give
\[\sumn_{n\le 3X} g_{z_\ell}(n) \ll_\delta X \Big(\frac{\log X}{\log z_\ell}\Big)^{80} \ll_\delta X\gamma^{-160\ell}.\]
If $1<z_\ell<2$, then $\log z_\ell < \log2$, and hence $\log X \ll_\delta \gamma^{-2\ell}$. 
Combining this with \eqref{eq:gz}, again yields same estimate for $\sumn_{n \le 3X} g_{z_\ell}(n)$ for  $1<z_\ell<2$. 
Hence for all $1<z_\ell \le 3X$, 
\[\sumn_{n\le 3X}g_{z_\ell}(n)\ll_\delta X\gamma^{-160\ell}.\]
Using \eqref{eq:w-eightyone} and the definition of $g_z$,  we obtain
\[\sumn_{n \le 3X} w(n)^2 \ll_{\delta,s} X + \sumn_{\ell\ge0} 4^{-\ell} \sumn_{n\le 3X} g_{z_\ell}(n) \ll_{\delta,s} X\Big( 1 + \sumn_{\ell\ge0} 4^{-\ell} \gamma^{-160\ell}\Big).\]
Now choosing $\gamma$ accordingly (for example, $\gamma>2^{-1/1000}$ works), the above infinite series converges, and therefore, we have the required estimate.
\end{proof}

\subsection{Upper bound for sign-weighted \texorpdfstring{$w$}{w}}\label{se:upper-bound-signed}
For simplicity, we write
\begin{align}\label{eq:big-w}
    W(n)=\sgn(c_f(n))w(n),
\end{align}
with $\sgn(0)=0$. Then the identity $\mu^2(ab)=\mu^2(a)\mu^2(b)\ind_{(a,b)=1}$ and \eqref{eq:bilinear} give
\begin{equation}\label{eq:W} 
    W(n) = \sum_{\substack{ab=n\\a\le y,\ V(a)>0,~(ab,Q)=(a,b)=1}} \mu^2(a) \mu^2(b) \rho^+(b) \sum_{1 \le i \le s} \beta_i(a) \lambda_i(b). 
\end{equation}

Recall that the indicator function of the square-free integers can be written as $\mu^2(n) = \sum_{d^2|n} \mu(d)$. For $D\ge 1$, define 
\[\mu^2_{\leq D}(n) := \sumn_{\substack{d^2 |n,~ d \leq D }} \mu(d)  \qquad \mu^2_{> D}(n) := \sumn_{\substack{d^2 |n,~ d > D }} \mu(d).\] 
With this, we  introduce the truncated weights 
\begin{equation}\label{W-small}
    W_{D}(n):= \sum_{\substack{ab=n\\a\le y,\ V(a)>0,~(ab,Q)=(a,b)=1}} \mu^2(a) \mu^2_{\leq D}(b) \rho^+(b)\sum_i\beta_i(a)\lambda_i(b), 
\end{equation}
and 
\begin{equation}\label{W-large}
    W^{D}(n):= W(n) - W_D(n) = \sum_{\substack{ab=n\\a\le y,\ V(a)>0,~(ab,Q)=(a,b)=1}} \mu^2(a) \mu^2_{> D}(b) \rho^+(b)\sum_i\beta_i(a)\lambda_i(b).
\end{equation}

With the notation above, we have the following upper bounds.
\begin{lem}\label{lem:tail}
For every $\eps>0$, we have
\begin{equation}\label{eq:pointwise-W}
    |W(n)|+|W_D(n)|\ll_{\eps,f}X^\eps \qquad(n\le4X),
\end{equation}
and
\begin{equation}\label{eq:tail:L1}
    \sum_{n\le4X} |W^D(n)| \ll_{\eps,f} \frac{X^{1+\eps}}{D}.
\end{equation}
\end{lem}

\begin{proof}
We first prove the pointwise bounds. Recall that  $\rho^+(b) \le \tau(b)$, and by Cauchy-Schwarz we also have
\[|\sumn_i \beta_i(a)\lambda_i(b)|\le\sqrt{s}\tau(b).\] 
Applying triangle inequality in \eqref{eq:W} and noting that the sum over $a\mid n$ contains at most $\tau(n)$ terms, from $|\mu^2_{\leq D}(b)|\le\tau(b)$ it follows that
\[W(n) \ll\tau(n)^3, \qquad W_D(n)\ll \tau(n)^4.\]  
This proves \eqref{eq:pointwise-W}.

For \eqref{eq:tail:L1}, we use the inequality $|\mu^2_{> Y}(n)| \leq \sum_{\substack{d^2 |n \\ d > Y }} 1.$
Then, 
\[|W^D(n)|=|W(n)-W_D(n)| \leq \sqrt{s} \sumn_{\substack{ab=n,~a\leq y}} \tau(b)^2 \sumn_{\substack{d^2 | b,~ d > D }} 1\] 
after dropping all coprimality and support restrictions.
Therefore,
\begin{align}
\sum_{n \le 4X} |W^D(n)| &\ll_{\eps,f} X^{\frac{\eps}{2}}~\sum_{a\le y}~\sum_{b\le4X/a}~\sum_{\substack{t>D\\ t^2 \mid b}}1 ~\ll~ X^{\frac{\eps}{2}}~\sum_{a\le y}~\sum_{\substack{t>D \\ t^2 \leq \frac{4X}{a}}} \left[\frac{4X}{at^2} \right]  \nonumber\\
&\ll X^{1+\eps/2} ~ \sumn_{a\le y} \frac{1}{a} \sumn_{\substack{t>D}} \frac{1}{t^2} ~  \ll ~ \frac{ X^{1+\eps/2} \log(2y)}{D} .\label{eq:tail-count}
\end{align}
Now writing $\log(2y)\ll \log X \ll X^{\frac{\eps}{2}}$ completes the proof of the lemma. 
\end{proof}

We next record the shifted-convolution bound which is the mixed-form analogue of \cite[Lemma~6.1]{mat15}, and follows from the same argument.
\begin{lem}\label{lem:shift}
There exists a constant $\tau>0$ such that, uniformly for positive integers
$a,b,A,B\le X^\tau$ and integers $r$ with $0<|r|\le X^\tau$, one has
\begin{equation}\label{eq:shift}
    \sumn_{\substack{X\le aAm,\; bBn\le 2X\\aAm-bBn=r}}~ \lambda_i(Am) \overline{\lambda_j (Bn)} \ll_f X^{1-\tau},
\end{equation}
uniformly over every pair $f_i$ and $f_j$.
\end{lem}

\begin{proof}
    For $A\geq 1$, let $f_{i,A}$ denote the cusp form obtained from $f_i$ by retaining only the Fourier coefficients whose indices are divisible by $A$.
    Similarly, let $\widetilde{f_j}$ be the cusp form such that $a_{\widetilde{f}_j}(n) = \overline{\lambda_j(n)}n^{(\kappa-1)/2},$ and define $\widetilde f_{j,B}$ from $\widetilde {f_j}$ with the coefficients whose indices are divisible by $B$. 
    Then $f_{i,A}$ and $\widetilde f_{j,B}$ are cusp forms with level growing  at most polynomially in $A$ and $B$, respectively. 
    The argument of \cite[Lemma~6.1]{mat15} now applies without change, as the estimate of \cite[Theorem 1]{har03} is valid for cusp forms of arbitrary levels and nebentypus. 
    Moreover, as noted in the addendum to \cite[Section~1]{har03}, the implied constant depends at most polynomially on the levels and the weights of the two forms. Therefore, for some sufficiently small $\tau_{i,j}>0$, we have
    \[\sumn_{\substack{X\le aAm,\; bBn\le2X\\ aAm-bBn=r}} \lambda_i(Am) \overline{\lambda_j(Bn)} \ll_{f_i,f_j} X^{1-\tau_{i,j}}\]
    uniformly for $a,b,A,B,r\le X^{\tau_{i,j}}$. 
    Since there are only finitely many pairs $f_i$ and $f_j$, after decreasing the exponent if necessary, we have the assertion.
\end{proof}

Now we prove the main proposition of this section.
\begin{prop}\label{prop:correlation} 
There exist $\tau>0$ and an absolute $\theta>0$ depending only on $f$ such that, uniformly for every  $0<|r|< X^\tau$, 
\begin{equation}\label{eq:corr} 
    \sumn_{\substack{X\le n, n+r\le 2X}} W(n)\overline{W(n+r)} \ll_f X^{1-\theta}. 
\end{equation}
\end{prop}

\begin{proof}
Let $C(r) = \sum_{X\le n, n+r\le 2X} W(n) \conj{W(n+r)}$ and $C_D(r) = \sum_{X\le n, n+r\le 2X} W_D(n)\conj{W_D(n+r)}.$
Then  the identity 
\[W(n) \conj{W(n+r)} - W_D(n) \conj{W_D(n+r)} = W^D(n) \conj{W(n+r)} + W_D(n) \conj{W^D(n+r)}\]
 and Lemma \ref{lem:tail} give 
\begin{equation}\label{eq:cr-estimate}
   C(r) = C_D(r) + O_{f,\eps} \left( \tfrac{X^{1+\eps}}{D} \right). 
\end{equation}

From the right-hand side of \eqref{W-small}, for $n\in [X,2X]$, we can write $W_D(n)$  as 
\[\sum_{\substack{a\le y,\ V(a)>0\\ (a,Q)=1}} \mu^2(a) \sum_{d|a} \sum_{u|Q} \sum_{e\in \mathcal{D}^+} \sum_{t\leq D} \mu(d)\mu(u)\mu(e)\mu(t) \sum_{\substack{ab=n\\B|b}} \sumn_i \beta_i(a)\lambda_i(b)\]
where $B= \lcm(d,u,e,t^2)$. 
 Then $C_D(r)$ equals 
\begin{align}\label{eq:cdr}
&\sum_{\substack{a, a'\le y\\ V(a), V(a')>0\\(a,Q)=(a',Q)=1}} \mu^2(a) \mu^2(a') \sum_{i,j} \beta_i(a) \conj{\beta_j(a')} \nonumber\\
&\times \sum_{\substack{d|a\\d'|a'}}  \sum_{\substack{u,u'|Q}} \sum_{\substack{e, e'\in \mathcal{D}^+}} \sum_{\substack{t, t'\leq D}}  \mu(d) \mu(d') \mu(u) \mu(u') \mu(e) \mu(e') \mu(t) \mu(t') \times \sum_{\substack{ab, a'b' \in [X,2X]\\ a'b'-ab=r\\B|b, ~B'|b'}} \lambda_i(b) \conj{\lambda_j(b')},
\end{align}
with $B'=\lcm(d',u',e',t'^2).$ 
Writing $b=Bm$ and $b'=B'm'$, the innermost sum is
\[\sum_{\substack{ X\leq aBm,\,a'B'm'\leq 2X\\ a'B'm'-aBm=r}} \lambda_i(Bm)\overline{\lambda_j(B'm')}. \]
By Lemma \ref{lem:shift}, the above is $\ll X^{1-\tau}.$ 
Recall that  $|\mathcal{D}^+|\leq y$ and the $\beta$-factors have absolute value $\le 1$. 
Also, by interchanging the summations, it is easy to see that  $\sum_{\substack{a\leq y}} \sum_{\substack{d|a}} 1 \ll y\log(2y).$ 
Thus, taking $D=X^\delta$ and dropping all the remaining restrictions in \eqref{eq:cdr}, we have
\begin{equation}\label{eq:cdr-estimate}
    C_D(r)\ll  X^{1-\tau} y^2 D^2 (y\log(2y))^2 \ll X^{1-\tau + 6\delta +\eps}.
\end{equation}

In order to apply Lemma \ref{lem:shift}, we need $B,B' \le X^\tau$. 
Note that $ B\le du e t^2\le Qy^2D^2=QX^{4\delta}$ (and similarly for $B'$). 
Thus, for $\delta < \tau$ sufficiently small, we have $B, B'\le X^\tau$ for all sufficiently large $X$, and $a,a'\le y< X^\tau$.  
Combining \eqref{eq:cr-estimate} and \eqref{eq:cdr-estimate} gives
\begin{equation}\label{eq:two-errors} 
    \sumn_{\substack{X\leq n, n+r \leq 2X}} W(n) \conj{W(n+r)} \ll_{\eps,f}X^{1-\tau+6\delta+\eps}+X^{1-\delta+\eps}. 
\end{equation}
By choosing $\eps>0$ sufficiently small in terms of $\delta$, each exponent on the right-hand side of \eqref{eq:two-errors} is at most $1-\theta$ with $\theta<\tau$.
Hence
\[\sumn_{\substack{X\leq n, n+r \leq 2X}} W(n) \conj{W(n+r)} \ll_f X^{1-\theta},\]
which proves the proposition.
\end{proof}

\subsection{Short interval estimates}
First we show that the bounds proved in Section \ref{se:lower-bound-w'} give the required lower bound for $w'$ on a positive proportion of short intervals. Set
\begin{align}\label{eq:eh}
    E_h = \Big\{x \in[X,2X] : \sumn_{x< n \le x+h}w'(n) \ge \frac{c_0h}{4}\Big\},
\end{align} 
where $c_0$ is the constant appearing in \eqref{eq:main-first}.

\begin{prop}\label{pr:lowerbound}
Let $f^\circ$ be non-lacunary. Then for any fixed $h \geq 1$, one has $\meas(E_h)\gg X$.
\end{prop}

\begin{proof}
From Proposition \ref{pr:moments}, we can choose constant $B>0$ such that, for sufficiently large $X$, 
\[ 0\le w'(n)\le B\quad(x< n\le 3X). \] 

As the interval $(x,x+h]$ contains at most $h+1$ integers, one has 
\begin{equation}\label{eq:A-upper} 
    \sumn_{x< n \le x+h} w'(n) \le B(h+1) \le 2Bh.
\end{equation} 
Thus, for sufficiently large $X$, by \eqref{eq:main-first} we obtain
\begin{align}\label{eq:A-average}
\int_X^{2X} \sumn_{x< n \le x+h}w'(n) \,dx & \ge \sumn_{n\ge 1,~x< n\le 2X} w'(n) \int_{n-h}^{n} dx \ge h\sumn_{X+h<n\le2X}w'(n) \\
&= h \Big(\sumn_{x< n\le2X} w'(n) - \sumn_{x< n \le x+h} w'(n) \Big) \\
&\ge h\bigl(c_0X-B(h+1)\bigr)\ge\frac{c_0}{2}hX.
\end{align}
Combining \eqref{eq:A-upper} and \eqref{eq:A-average}, and splitting the integral over $[X,2X]$ into contributions from $E_h$ and its complement, we obtain
\[\frac{c_0}{2}hX\le\int_X^{2X}\sumn_{x< n \le x+h}w'(n)\,dx\le\frac{c_0}{4}hX+2Bh\meas(E_h).\]
Rearranging yields
\begin{equation}\label{eq:measure-lower-bound}
    \meas(E_h)\ge\frac{c_0}{8B}X \gg X. \qedhere
\end{equation} 
\end{proof}

An application of Proposition \ref{pr:upperbound2}, combined with Proposition \ref{prop:correlation}, yields the following result. 
For this, let 
\begin{align}\label{eq:bk}
    B_K = \Big\{ x \in [X,2X] : \Big|\sumn_{x < n \le x+h} \sgn(c_f(n))w(n)  \Big| > K\sqrt {h}\Big\}.
\end{align}

\begin{prop} \label{pr:upperbound}
For every $K>0$ and fixed $h\ge 1$, one has $\meas(B_K) \ll_f \frac{X}{K^2}.$
\end{prop}
\begin{proof}
    We follow the proof of \cite[Lemma~6.2]{mat15}. Expanding the square, writing $n_2 - n_1=r$ in $\sum_{x< n_1 \le x+h} \sum_{x< n_2 \le x+h}$, and interchanging the integral and the finite sums (cf. proof of \cite[loc. cit]{mat15}) gives
    \begin{align}
        \int_X^{2X} \Big| \sumn_{x< n \le x+h}W(n) \Big|^2\,dx = \sumn_{X\le n\le2X}~ \sumn_{|r|\le h} (h-|r|)W(n)W(n+r) + E,
    \end{align}
where the error term $E$ is supported on pairs for which both variables lie within distance $h$ of one of the endpoints $X$ or $2X$. 
More precisely,
\[|E| \ll h \sumn_{\substack{|n_i-X|\le h}} |W(n_1)W(n_2)| + h \sumn_{\substack{|n_i-2X|\le h}} |W(n_1) W(n_2)|.\]
By the pointwise bound $|W(n)| \ll_{\varepsilon,f} X^\varepsilon$ (see the proof of Lemma \ref{lem:tail}) for $n\leq 3X$, there are $O(h^2)$ pairs in each of the two sums, and hence
\begin{equation}\label{eq:endpoint-error}
    E \ll_{\varepsilon,f} h^3 X^{2\varepsilon}.
\end{equation}
We next separate the diagonal $r=0$. Since $|W(n)| = w(n)$, we have
\begin{align}\label{eq:main-term-1}
    h \sumn_{X \le n\le2X} W(n)^2 \le h \sumn_{n \le 3X} w(n)^2.
\end{align}

For the off-diagonal part, note that the terms for which one of $n,n+r$ falls outside $[X,2X]$ lie within distance $h$ of an endpoint and are already covered by an error of the same type as \eqref{eq:endpoint-error}. 
Thus, using \eqref{eq:corr},
\begin{align}\label{eq:main-term-2}
    \sum_{0<|r|\le h} (h-|r|)~ \Big| \sum_{\substack{X\le n \le 2X\\ X \le n+r \le 2X}} W(n) W(n+r) \Big| &\ll_f X^{1-\theta} \sum_{0<|r|\le h}(h-|r|) \ll_f h^2 X^{1-\theta}. 
\end{align}
Combining \eqref{eq:endpoint-error}, \eqref{eq:main-term-1} and \eqref{eq:main-term-2}, we obtain
\begin{align}\label{eq:shifted-big-w}
    \int_X^{2X} \Big| \sumn_{x< n \le x+h} W(n) \Big|^2\,dx \ll h^3X^{2\varepsilon} + h\sumn_{n\le 3X}w(n)^2 + h^2X^{1-\theta}.
\end{align}
Since $\theta>0$ and $h$ is fixed, from Proposition \ref{pr:upperbound2} it follows that the right-hand side is $O_{f}( hX).$ 
Now applying Chebyshev's inequality, we immediately have
\begin{align}
    C_f hX \ge \int_X^{2X} \Big| \sumn_{x< n \le x+h} \sgn(c_f(n)) w(n) \Big|^2\,dx \ge K^2 h \meas(B_K),
\end{align} 
for some constant $C_f>0$, and the assertion follows.
\end{proof}
\subsection{Proof of \protect\eqref{asymp-nonCM} --  conclusion}\label{se:proof-noncm}
In this section we finally prove $\displaystyle \mathscr{S}^{\flat}_f(X) \asymp \mathscr{N}_f^\flat(X)$. Let us define the set
\begin{align} \label{}
    D(x,h):= \{ x \in [X,2X] \colon \Big | \sumn_{x \leq n \leq x + h} \text{sgn}(c_f(n)) w(n) \Big | < \sumn_{x \leq n \leq x + h} w'(n) \}.
\end{align}
Recall the set $B_K$ from \eqref{eq:bk}. Then from Proposition~\ref{pr:upperbound}, we get that
$ \displaystyle \meas(B_K) \le C_f X/K^2$. 
Fix $K>0$ so large that $C_f K^{-2}<\frac{c_0}{16B}$.

Further recall the set $E_h$ from \eqref{eq:eh}. From Proposition~\ref{pr:lowerbound}, fix an integer $h$ so large that $K \sqrt{ h}<\frac{c_0}{4}h $. 
Both the choices of $K,h$ are independent of $X$. 

Note that by our choices of parameters, we have
\begin{align}
   E_h \setminus B_K \subseteq D(x,h).
\end{align}
Therefore, by \eqref{eq:measure-lower-bound} and Proposition \ref{pr:upperbound}, $\mathrm{meas}(D(x,h))$ is at least $\frac{c_0 X}{16B}\gg X$.
For every such $x$, a sign change of $c_f(n)$ must occur in the interval $[x, x+h]$, since $w(n) \geq w'(n) \geq 0$ for every $n$. 
Hence, from Lemma \ref{le:sign-change-count}, we have  $\mathscr{S}_f^\flat(X) \gg X$ for sufficiently large $X$. 
Since  in the current setting $\mathscr{N}_f(X)\asymp X$, we have 
\[\mathscr{N}_f(X) \gg \mathscr{N}_f^\flat(X)  \gg \mathscr{S}_f^\flat(X) \gg \mathscr{N}_f(X).\] 

\section{\texorpdfstring{Notation and setting in the case $\kappa \in \Z$: $ f^\circ$}{f ◦ } is lacunary}\label{se:all-cm}

When $f^\circ$ is lacunary, the situation is more subtle than what we have considered so far because this includes both CM-forms and weight 1 cusp forms (these are always lacunary). 
We need to set-up further notation in order to proceed.
For details on weight 1 cusp forms, we refer the reader to \cite{DS74, Serre81}, and \cite{Ribet1977, miyake} for CM forms.

Note that $f^\circ \neq 0$ since $f$ satisfies \ref{tag:1}. 
For $Q=2N$ and each $j$ appearing in \eqref{eq:newform-decomp-simplified-coeffs} such that $\alpha_{j} = \alpha_{j,1} \neq 0$, we define 
\begin{gather}\label{eq:defn-aj}
     Z_j:= \{p\nmid Q:\lambda_j(p) = 0\}, \quad \mathcal{A}_j = \{n\ge 1: \mu^2(n)=1, (n,Q)=1, \lambda_j(n)\neq 0\}, 
     \\
     \mathcal{A}_j(t) = \{n\leq t: n\in \mathcal{A}_j \}.
\end{gather}

\subsection{\texorpdfstring{$\kappa=1$}{κ=1}}
Let $f\in S_1(\Gamma_1(N))$ be a cusp form with Fourier expansion \eqref{expansion} such that $a_f(n)$ are real for all $n$. 
Further, assume that $f$ satisfies \ref{tag:1}. 
By Deligne-Serre \cite[Theorem 4.1]{DS74}, for each $j$ appearing in  \eqref{eq:newform-decomp-simplified}, one can associate to $f_j$ an irreducible odd representation $\rho_j :\mathrm{Gal}(\overline{\Q}/\Q) \longrightarrow \mathrm{GL}_2(\C)$ with finite image $G_j$, satisfying $\lambda_j(p)=\operatorname{tr}\rho_j(\operatorname{Frob}_p)$. 
Define
\begin{align}\label{eq:notation-wt1}
     C_j := \{g\in G_j: \mathrm{tr}(g)=0\}, \quad \theta_j := {|C_j|}/{|G_j|}, \quad \theta_f := \min_{\alpha_{j,d}\neq 0} \theta_j, \quad \text{and} \quad \theta_{f^\circ} = \theta_{\circ} :=\min_{\alpha_{j,1}\neq 0}\theta_j.
\end{align}
Then, $0<\theta_j <1$ for each $j$. Note that if $a_f(n)\neq0$, then for some $(j,d)$ with $\alpha_{j,d}\neq0$, we have $d\mid n$ and $\lambda_j(n/d)\neq0$. 
Hence
    \begin{align}\label{eq:nf-count-wt1}
        \mathscr{N}_f (X) \le \sumn_{\alpha_{j,d} \neq 0} \#\{m \le X/d : \lambda_j(m) \neq 0\} \ll   X/(\log X)^{\theta_f}
    \end{align}  
    by \cite[\S 7.7, Th\'eor\`eme 18]{Serre81}.
    Similarly, we have 
    \begin{align}\label{eq:nf-flat-count-wt1}
        \mathscr{N}_f^\flat(X) \ll X(\log X)^{-\theta_{\circ}}.
    \end{align}

Note that $\mathrm{Frob}_p\in C_j$ for all $p\in Z_j$, and since $C_j$ is invariant under conjugation, by an effective form of Chebotarev (see \cite[\S, 2.2, Th\'eor\`eme 2]{Serre81}),
\begin{align}\label{eq:density-zj}
    \#\{p\le x: \lambda_j(p)=0\} = \theta_j {x}/{\log x} + O(x(\log x)^{-2}).
\end{align}
Also note that, from multiplicativity and \eqref{eq:defn-aj}, it follows that 
\begin{align}\label{eq:ajt-wt1}
    \mathcal{A}_j = \{n\ge 1: \mu^2(n)=1, (n,Q)=1, p\mid n \Longrightarrow p \notin Z_j\}.
\end{align}
Since $\ind_{\mathcal A_j}$ is a non-negative multiplicative function, by \cite[\S 6.1, Th\'eor\`eme 13]{Serre81} we have 
    \begin{align}\label{eq:ajt-wt1-bound}
        \mathcal{A}_j(t) \asymp_f {t}{(\log t)^{-\theta_j}}.
    \end{align}

\subsection{\texorpdfstring{$\kappa > 1$}{κ > 1}}
For an imaginary quadratic field $K$ define 
\[\mathscr{A}_K := \left\{ n\geq 1 ~\colon~ \mu^2(n)=1,\ (n,Q)=1,\ \text{and } p\mid n \Longrightarrow p \text{ splits in } K \right\}, \quad \mathscr A_K(t):=\#\{n\leq t~\mid~n \in \mathcal{A}_K\}.\]

We assume that each $f_j$ occurring in \eqref{eq:newform-decomp-simplified-coeffs} has CM by an imaginary quadratic field $K_j$.
Then, we have $\mathcal{A}_j = \mathscr{A}_{K_j}$ (follows from Lemma~\ref{lem:psplit}) and $Z_j=\{p\nmid Q : p \text{ is inert in } K_j\}$.
Moreover, by \cite[Th\'eor\`eme 13]{Serre81}, we have 
\begin{align}\label{eq:akt-cm-bound}
    \mathcal{A}_j(t) = \mathscr A_{K_j}(t) \asymp t(\log t)^{-1/2}.
\end{align}

\begin{lem} \label{lem:psplit}
    Let $f \in S_k(N)$ be a  newform of  weight $k\geq 2$.  Suppose that $f$ has CM by the imaginary quadratic field $K$. Let $p\nmid N$ be a rational prime which splits in $K$. Then $\lambda_f(p) \neq 0$.
\end{lem}

\begin{proof}
By the theory of CM newforms (see e.g. \cite{Ribet1977}, \cite{miyake}), there is a primitive algebraic Hecke character $\displaystyle \psi:I_K(\mathfrak f)\longrightarrow \mathbb C^{\times}$ of conductor $\mathfrak f$ and infinity type $(k-1,0)$ whose theta series is $f$. 
Here $I_K(\mathfrak f)$ denotes the group of fractional ideals of $K$ relatively prime to $\mathfrak f$. 
With one standard choice of convention, the infinity-type condition says that
\begin{equation}\label{eq:principal-character}
    \psi((\alpha))=\alpha^{k-1}  \qquad \text{if }\alpha\equiv 1\pmod{\mathfrak f}.
\end{equation}
If the conjugate convention is used, the right-hand side is $\overline\alpha^{\,k-1}$; the proof below is unchanged.

The Fourier expansion of $f$ is thus $\displaystyle \sum_{\mathfrak a\subseteq \mathcal O_K\ \ (\mathfrak a,\mathfrak f)=1 } \psi(\mathfrak a)q^{N(\mathfrak a)}$, where $N(\cdot) = N_{K|\Q}(\cdot)$ denotes the field norm.
Suppose that $p\nmid N$ splits as $\displaystyle p  \mathcal O_K=\mathfrak p\overline{\mathfrak p}$.
The primes $\mathfrak p$ and $\overline{\mathfrak p}$ are distinct and are relatively prime to $\mathfrak f$. 
They are exactly the integral ideals of $K$ having norm $p$. 
We have the relation,
\begin{equation}\label{eq:split-coefficient}
    a_f(p)=\psi(\mathfrak p) + \psi(\overline{\mathfrak p}).
\end{equation}

Assume, for a contradiction, that $a_f(p)=0$. 
Then
\begin{equation}\label{eq:ratio-minus-one}
    \frac{\psi(\mathfrak p)} {\psi(\overline{\mathfrak p})}=-1.
\end{equation}
Consider the fractional ideal $\mathfrak b=\mathfrak p\overline{\mathfrak p}^{\,-1} \in I_K(\mathfrak f)$.
The ray class group modulo $\mathfrak f$ is finite.  Hence there are an integer $r\geq 1$ and an element $\alpha\in K^{\times}$ such that
\begin{equation}\label{eq:ray-principalization}
     \mathfrak b^r=(\alpha), \qquad  \alpha\equiv 1\pmod{\mathfrak f}.
\end{equation}
Using multiplicativity of $\psi$, equations \eqref{eq:principal-character} and \eqref{eq:ray-principalization} give
\begin{align*}
 \left(   \frac{\psi(\mathfrak p)} {\psi(\overline{\mathfrak p})} \right)^r =  \psi(\mathfrak b^r)  =\psi((\alpha))  =\alpha^{k-1}.
\end{align*}
By \eqref{eq:ratio-minus-one}, we obtain
\[\alpha^{k-1}=(-1)^r,  \qquad\text{and therefore}\qquad \alpha^{2(k-1)}=1.\]
Since $k-1>0$, the element $\alpha$ is a root of unity.  
Consequently $(\alpha)= \mathcal O_K$. On the other hand,
\[ (\alpha)=\mathfrak b^r  =\mathfrak p^r\overline{\mathfrak p}^{\,-r} \neq  \mathcal O_K,\]
because its $\mathfrak p$-valuation is $r$ and its $\overline{\mathfrak p}$-valuation is $-r$.  
This contradiction proves that $a_f(p)\neq 0$ whenever $p\nmid N$ splits in $K$.
\end{proof}

\section{Proof of Theorem \ref{thma} when \texorpdfstring{ $\kappa \in \Z$: $ f^\circ$}{f ◦ } is lacunary}\label{se:all-cm-11}
In this section we prove that when $f^\circ$ is lacunary, 
\begin{align} \label{asymp-allCM}
\mathscr{S}_f^\flat(X) \asymp \mathscr{N}_f^\flat(X) \asymp
X(\log X)^{-\theta_{\circ}},
\end{align}
where $\theta_{\circ}$ is defined by \eqref{eq:theta-f-circ}.
We will borrow the notation, and use the information gathered in Section \ref{se:positive} to this end. The proof is again by comparing the sieve weights $w,w'$ along the non-zero Fourier coefficients $c_f(n)$.

\subsection{Lower bound for \texorpdfstring{$w'$}{w′}}
Throughout the rest of the section, $f\in S_\kappa(\Gamma_1(N))$, $\kappa\ge 1$. Recall that 
\begin{align}\label{eq:w'-redefine}
w'(n) = \mu^2(n) \ind_{(n,Q)=1} \sum_{\substack{ab=n\\a\le y,\ V(a)>0\\P^-(b)>y}} \frac{|c_f(n)|}{V(a)}; \quad H_V(y)= \sum_{\substack{a\leq y,~  (a,Q)=1\\ V(a)>0}} \frac{\mu^2(a)}{a},
\end{align} 
and from \eqref{eq:defn-aj}, \eqref{eq:ajt-wt1-bound} and \eqref{eq:akt-cm-bound}, 
\begin{align}\label{eq:support-cm}
    \left\{ n\geq 1 ~\mid~ \mu^2(n)=1,\ (n,Q)=1,\ c_f(n)\neq 0 \right\} \subseteq \bigcup\nolimits_{j}\mathcal{A}_{j}, \quad \mathcal{A}_{j}(t) \ll t (\log t)^{-\theta_{\circ}},
\end{align}
where $\theta_{\circ}$ is given by \eqref{eq:notation-wt1} if $\kappa =1$, and $\frac{1}{2}$ otherwise. i.e.,
\begin{align}\label{eq:theta-f-circ}
    \theta_{\circ} = \begin{cases}
        \min_{\alpha_{j,1}\neq 0}\theta_j, \quad &\text{ if $\kappa =1$,}\\
        \frac{1}{2}, \quad &\text{ if $\kappa\ge 2$.}
    \end{cases}
\end{align}
\begin{lem}\label{le:cm}
If $f^\circ$ is lacunary, then $H_V(y) \gg (\log y)^{1-\theta_{\circ}}.$
\end{lem}
\begin{proof}
Since $V(a)\ge |\alpha_i \lambda_i(a)|$ for any $i$, we have the trivial inequality
\[\sumn_{\substack{a\leq y,~  (a,Q)=1,~ V(a)>0}}\frac{\mu^2(a)}{a} \geq \sumn_{\substack{a\leq y,~ a \in \mathcal A_i}} \frac{1}{a}. \] 

On the other hand, from the asymptotics in \eqref{eq:ajt-wt1-bound} and \eqref{eq:akt-cm-bound}, we have
\begin{equation}\label{eq:harmonic-cm}
\sumn_{\substack{a\le t,~ a\in\mathcal A_i}}\frac{1}{a} = \frac{\mathcal A_i(t)}{t}+\int_1^t\frac{\mathcal A_i(u)}{u^2} du \asymp (\log t)^{1-\theta_{i}}.
\end{equation}
By applying this to $f_i$ the claim now follows. 
\end{proof}

\begin{prop}\label{pr:moments-cm}
    Suppose that $f^\circ$ is lacunary. Then
    \begin{equation}\label{eq:main-first-cm}
        \sumn_{X < n\le 2X} w'(n) \gg_f X (\log X)^{-\theta_\circ}.
    \end{equation}
\end{prop}

\begin{proof}
    The assertion follows easily from Lemma \ref{le:first-moment-w'-preli} and Lemma \ref{le:cm}.
\end{proof}

\subsection{Upper bound for \texorpdfstring{$w$}{w}}
The next proposition controls the $L^2$-norm of $w$ in the lacunary setting.
\begin{prop}\label{pr:upperbound2-cm}
    Let $f^\circ$ be lacunary. Then for sufficiently large $X,$ we have $$\sumn_{n\le 3X} w(n)^2 \ll_f X(\log X)^{-\theta_\circ}.$$
\end{prop}
\begin{proof}
    We have already obtained \eqref{eq:w-eightyone}, which does not depend on the CM-type of the forms $f_j$. 
    Since $c_f(n)=0$ for $n\notin \bigcup_{j}\mathcal{A}_{j}$, from \eqref{eq:w-eightyone} we can write 
    \begin{equation}\label{eq:supportmajorant-cm}
    \sumn_{n\le 3X} w(n)^2 \ll_{\delta,s} \sumn_i \sumn_{n\le 3X} \ind_{\mathcal A_i}(n) + \sumn_{\ell\ge0}4^{-\ell} \sumn_i \sumn_{n\le 3X} \ind_{\mathcal A_{i}}(n) 81^{\omega_{>z_\ell}(n)}.
\end{equation} 
Similar to the proof of Proposition \ref{pr:upperbound2}, for $z>1$, we define the non-negative multiplicative function
\[g_{i,z}(n) := \ind_{\mathcal A_{i}}(n)\cdot 81^{\omega_{>z}(n)}.\] 
Then, for a prime $q$, we have
\[g_{i, z}(p)= \begin{cases}
0, & \text{if } p \in Z_i,\\
1, & \text{if } p\nmid Q, p\notin Z_i, \text{ and } p\le z,\\
81, & \text{if } p\nmid Q, p\notin Z_i, \text{ and } p> z.
\end{cases}\]
Note that $g_{i,z} (p^\nu)=0$ for any $\nu\ge 2.$
Thus, from Lemma \ref{le:MVmultfunct}, we have
\begin{align}\label{eq:gz-cm}
    \sum_{n \le 3X} g_{i,z}(n) &\ll X \prod\nolimits_{\substack{q\leq z\\ q  \in Z_i}} (1 - \tfrac 1q) \prod\nolimits_{z\leq q\le 3X} \left(1-\tfrac1q + \tfrac{g_{f_i,z}(q)} q \right) \nonumber \\
    &\ll \frac{X}{(\max(1,\log z))^{\theta_i}} \left(\frac{\log(6X)}{\max(1,\log z)}\right)^{80} \ll \frac{X}{(\log X)^{\theta_i}} \left(\frac{\log(6X)}{\max(1,\log z)}\right)^{80},
\end{align} 
where we have used $\log(1+t)\le t$ along with 
\begin{align}\label{eq:merten-type}
    \sumn_{\substack{p<z,~ p\in Z_j}} p^{-1} = \theta_j \log\log z+O_{f_j}(1).
\end{align} 
For $\kappa \ge 2$, a proof can be found in \cite[Proposition~3.3.2]{pau23}. 
The case $\kappa = 1$  follows from \eqref{eq:density-zj} by partial summation.

Thus, for each fixed $\ell$, the second term in \eqref{eq:supportmajorant-cm} contributes
\begin{align}
    \sumn_i \sumn_{n\le 3X} g_{i, z_\ell}(n) \ll X(\log X)^{-\theta_i} \gamma^{-162\ell}.
\end{align}
Indeed, $\log z_\ell = (1-\gamma^2)\gamma^{2\ell}\log y$, so taking $z=z_\ell$ in \eqref{eq:gz-cm} and putting $y=X^\delta$ gives the claim. 
Also, the contribution from the first term is $\sum_i \sum_{n\le 3X} \ind_{\mathcal A_i}(n) \ll X(\log X)^{-\theta_i}$. 
Suitably choosing$\gamma$, the assertion now follows from \eqref{eq:supportmajorant-cm}.
\end{proof}
\subsection{Short-interval estimate}
To prove the analogue of Proposition \ref{pr:lowerbound} in the present case, we will need an auxiliary result. For an integer $d>0$, put
\begin{equation}\label{nfxd}
        N_j (X,d)    :=    \#\{n ~|~ X < n \leq 3X,\ n,\ n+d\in \mathcal{A}_{j}\}.
    \end{equation} 
We note that the following holds for any imaginary quadratic field $K$ and $\mathcal{A}_{K}$.
\begin{lem}\label{le:sieve-cm}
    We have
    \[N_j (X,d) \ll X(\log X)^{-2\theta_j} \prod\nolimits_{p|d} \Big( 1+ \tfrac{4}{p}\Big).\]
\end{lem}
\begin{proof}
    For $z\ge 3$, let 
    \begin{align}
        S_j(X,d;z) = \#\left\{ X < n\leq 3X:  \bigl(n(n+d),P_j(z)\bigr)=1 \right\},
    \end{align}
    where $P_j(z) = \prod\nolimits_{\substack{p\leq z\\ p\in Z_j}}p$.
    Then we have
\begin{equation}\label{eq:nfxdleqsxdz}
    N_j (X,d)\leq S_j(X,d;z).
\end{equation}
In what follows, we remove the subscript $j$ for simplicity. 
For any $p\in Z$, define $\Omega_d(p)=\{0,-d\}\subseteq \mathbb Z/p\mathbb Z$ to be the residue classes avoided by any $n$ counted in $N (X,d)$ and 
\[\nu_d(p) := |\Omega_d(p)| = \begin{cases}
2, & p\nmid d,\\
1, & p\mid d.
\end{cases}\]
For $q\mid P(z)$, define $\nu_d(q) = \prod_{p\mid q}\nu_d(p),$  $\nu_d(1)=1$. 
By the Chinese remainder theorem, each choice of one forbidden class at every prime $p\mid q$ determines exactly one residue class modulo $q$. 
Therefore, there are $\prod_{p\mid q}\nu_d(p)=\nu_d(q)$ forbidden residue classes modulo $q$, and summing over these $\nu_d(q)$ residue classes gives
\begin{align}
    \#\left\{ X < n\leq 3X: q\mid n(n+d) \right\} = \frac{2X\nu_d(q)}{q}+r_q, \qquad |r_q|<\nu_d(q) < \tau(q).
\end{align}

Define $g_d(p) = \nu_d(p)/p$ or $0$ according as $p\in Z$ or $p\notin Z$. 
Using $\log(1-t)\le -t$ for $0\le t <1$ , we have 
\begin{align}
    -\sumn_{w\le p <z}  \log(1- g_d(p)) \le  2\sumn_{w\le p <z}  p^{-1} +O(1),
\end{align}
where the implied constant is absolute. 
Thus, Mertens' estimate give 
\begin{align}\label{eq:5-38}
    \prod\nolimits_{w\le p <z} (1- g_d(p))^{-1} \le K (\log z/\log w)^2, \quad K>1.
\end{align}
Now, we apply \cite[Corollary 10]{fri10} to our situation with $\kappa =2$ to obtain
\begin{align}\label{eq:sxdz-leq-product-1}
    S(X,d;z)  \ll X \prod\nolimits_{\substack{p <z\\p\in Z}} (1- g_d(p)).
\end{align}
By definition of $g_d(p)$, the product on the right-hand side can be written as
\begin{align} \label{eq:log-ineq-cm}
    \prod_{\substack{p\mid P(z),~ p\nmid d}} \Big(1-\tfrac{2}{p}\Big) \prod_{\substack{p\mid (d,P(z))}} \Big(1-\tfrac 1p\Big) = \prod_{\substack{p\mid P(z)}} \Big(1-\tfrac{2}{p}\Big) \prod_{\substack{p\mid (d,P(z))}} \frac{p-1}{p-2}.
\end{align}

Using $\log(1-t)\leq -t$ for $0<t<1$ along with Mertens' estimate again, from \eqref{eq:merten-type} and \eqref{eq:log-ineq-cm}, we finally obtain 
\[\prod\nolimits_{\substack{p <z\\p\in Z}} (1- g_d(p))\ll (\log z)^{-2\theta} \prod\nolimits_{\substack{p\mid d}} \left(1+\frac4p\right).\]
Putting this back in \eqref{eq:sxdz-leq-product-1} with $D=X^{1/10}$ and $z=X^{1/100}$, we have the assertion of the lemma.
\end{proof}

Now we prove the following short-interval estimate for $w'$ in the lacunary case.
\begin{prop}\label{pr:lowerbound-cm}
Suppose that $f^\circ$ is lacunary. Then, for any $h \geq 1$ and sufficiently large $X$, there is a positive proportion of $x \in [X, 2X]$ such that
\begin{equation}\label{eq:shortlowbound-cm}
    \sumn_{x \leq n \leq x + h(\log X)^{\theta_\circ}} w'(n) \gg h.
\end{equation} holds.
\end{prop}
\begin{proof}
    From \eqref{eq:w'-redefine} it easily follows  that $w'(n) \leq 2^{\omega(b)} \sqrt{s} \cdot \ind_{\cup_i \mathcal{A}_i} ~(n) \leq 2^{\omega(b)}\sqrt{s} \sumn_i~ \ind_{{ \mathcal{A}_i}} (n),$ where $s$ is the number of non-zero $\alpha_j$ in \eqref{eq:newform-decomp-simplified-coeffs}.
    Let  
    \[H=h(\log X)^{\theta_\circ}, \quad \text{ and } \quad \mathscr{C}_i(x) = \#\{x\leq n \leq x+H ~|~ n\in \mathcal{A}_i\}.\] 
    Then, from the support condition \eqref{eq:support-cm} we have $\sum_{x\leq n\leq x+H} w'(n) \leq 2^{\omega(b)} \sqrt{s} \cdot \sumn_i \mathscr{C}_i (x)$, and consequently, 
    \begin{equation}\label{eq:secondemoment-shortsum-w'-cm}
        \int_X^{2X} \Big( \sumn_{x\leq n\leq x+H} w'(n) \Big)^2 dx \ll_{f,\delta} \sumn_i~ \int_X^{2X} \mathscr{C}_i^2(x)dx
    \end{equation} 
    by Cauchy-Schwarz.
    Note that the estimate on the right-hand side is now purely a counting problem with no dependence on $V$ or $y.$ 
    Now expanding the square gives 
    \begin{align}
        \mathscr{C}_i^2(x) &= \sumn_{x\leq m, n\leq x+H} \ind_{\mathcal{A}_i} (n) \ind_{\mathcal{A}_i} (m).
    \end{align}
    For $m=n,$ the condition $x< n\leq x+H$ is equivalent to $n-H\leq x< n$.
    Therefore the diagonal contribution to the  integral $\int_X^{2X} \mathscr{C}_i^2(x) dx$ is 
    \begin{align}\label{eq:c-square-first-part}
        \le H \mathcal A_i (3X) \ll {HX}(\log X)^{-\theta_i}.
    \end{align}     
    For the off-diagonal part, suppose that $m=n+d$ with $d>0$. 
    Such a pair contributes only when $d<H$, and thus from the above diagonal estimate we obtain 
    \begin{align}\label{eq:int-cfx}
          \int_X^{2X}\mathscr{C}_i (x)^2\,dx   & \leq    H \mathcal A_i(3X)    +    2\sum_{1\leq d<H}(H-d)N_i (X,d). 
    \end{align}
    From Lemma \ref{le:sieve-cm}, we infer that \begin{align*}
        \sum_{1\leq d <H} N_i (X,d) &\ll \frac{X}{(\log X)^{2\theta_i}} \sum_{1\leq d <H}  \prod_{\substack{p|d\\ p>2}} \Big( 1+ \frac{4}{p}\Big) =  \frac{X}{(\log X)^{2\theta_i}} \sum_{\substack{a\leq H\\ a \text{ odd}}} \frac{\mu^2(a) 4^{\omega(a)}}{a} \lfloor H/a\rfloor\\
        &\leq \frac{HX}{(\log X)^{2\theta_i}} \sum_{\substack{a\geq 1}} \frac{\mu^2(a) 4^{\omega(a)}}{a^2} = \frac{HX}{(\log X)^{2\theta_i}}\prod_{p} (1+4/p^2) \ll \frac{HX}{(\log X)^{2\theta_i}}.
    \end{align*} 
    Combining  the diagonal and off-diagonal contribution, the left-hand side of \eqref{eq:int-cfx} is $\ll XH^2(\log X)^{-2\theta_i}.$
    Since $\theta_i \ge \theta_0$, from \eqref{eq:secondemoment-shortsum-w'-cm} we  obtain
    \begin{equation}\label{eq:second-moment-w'}
        \int_X^{2X} \Big( \sumn_{x\leq n\leq x+H} w'(n) \Big)^2 dx \leq C_2 Xh^2.
    \end{equation} 
    for $h\geq1$, where $C_2>0$ is independent of $h$.

Let $X$ be sufficiently large. Then $H<X/4$, and for every integer $n$ with $5X/4 < n \le 7X/4$, the interval $[n-H, n)$ is contained in $[X,2X]$. 
Since $w'(n)\geq0$, we have 
    \begin{align}\label{eq:first-moment-cm}
        \int_X^{2X} \sumn_{x\leq n\leq x+H} w'(n)  dx \ge  H \sumn_{5X/4<n\leq7X/4}w'(n) \ge c_2 hX
    \end{align}
 for  a fixed constant $c_2>0$, where the last inequality follows  from Lemma \ref{le:first-moment-w'-preli} and Lemma \ref{le:cm}. 
 Let 
\begin{align}\label{eq:eh-cm}
    E_h' := \Big\{x\in[X,2X]~:~  \sumn_{x\leq n\leq x+H} w'(n) \geq \frac{c_2}{2} h\Big\}.
\end{align} 
Then splitting the integral  in \eqref{eq:first-moment-cm} into contributions from $E_h'$ and its complement, we have 
\[\int_{E_h'}  \sumn_{x\leq n\leq x+H} w'(n)\,dx \ge \frac{c_2}{2}hX.\]
Now applying Cauchy-Schwarz and \eqref{eq:second-moment-w'}, we can easily show that 
\[C_2h^2 X \meas(E_h')  \geq \frac{c_2^2 h^2 X^2}{4}.\] 
Thus $\meas(E_h')\gg X$, as desired.
\end{proof}

\subsection{Conclusion of the proof of \protect\eqref{asymp-allCM} -- Theorem~\ref{thma} for integral weights}
    Since the results of Section \ref{se:upper-bound-signed} remain valid regardless of whether all $f_j$ in \eqref{eq:newform-decomp-simplified-coeffs} are lacunary, from \eqref{eq:shifted-big-w} we obtain
    \[\int_X^{2X} \big| \sumn_{x< n \le x+ H} W(n) \big|^2\,dx \ll H^3 X^{2\eps} + H\sumn_{n\leq 3X} w(n)^2 + H^2 X^{1-\theta} \ll_f hX,\]
    where $H=h(\log X)^{\theta_\circ}.$ 
    For any $K>0$ and fixed $h\ge 1$, define
    \[B_K':= \Big\{ x\in [X,2X]~:~ \Big| \sumn_{x< n\leq x+H} W(n) \Big| > K\sqrt{h} \Big\}.\] 
    Now Chebyshev's inequality and the above estimate give 
    \[\mathrm{meas}(B_K') \ll_f \frac{X}{K^2}.\] 
    
    On the other hand, Proposition \ref{pr:lowerbound-cm} gives a set $E_h'$ with $$\meas(E_h') \gg_f X$$ such that \eqref{eq:shortlowbound-cm} holds for $x\in E_h'$. 
    Now proceeding exactly as in the proof of \eqref{asymp-allCM} in Section \ref{se:one-non-cm}, we obtain that a sign change must occur in the interval $(x,x+H]$ for each $x$ in a set of measure $\gg X.$ 
    Therefore, for sufficiently large $X$, from Lemma \ref{le:sign-change-count}, we must have 
    $$\mathscr S_f^{\flat}(X) \gg X/H \gg X(\log X)^{-\theta_\circ}.$$  
    On the other hand, since all the newforms in \eqref{eq:newform-decomp-simplified-coeffs} are lacunary, from \eqref{eq:serre-estimates}, \eqref{eq:theta-f-circ}, and \eqref{eq:nf-flat-count-wt1} we have $\mathscr N_f^{\flat}(X) \asymp  X(\log X)^{-\theta_\circ}$. 
    This completes the proof. \qed


\section{Proof of Theorem \ref{thma}} \label{sec:thma-end}
The different cases of Theorem \ref{thma} have already been established separately in Sections \ref{se:half-int}, \ref{se:one-non-cm}, and \ref{se:all-cm-11}. 
For the sake of completeness, we explicitly collect these results together here for the convenience of the reader. 

The half-integral weight case follows from Section \ref{se:half-int}.
For integral weights, $\mathscr{N}^\flat_f(X)$, and therefore $\mathscr{S}^\flat_f(X)$ only depend on the lacunarity of $f^\circ$. 
The case where $f^\circ$ is non-lacunary is treated in Section~\ref{se:one-non-cm} and the case where it is lacunary is treated in Section~\ref{se:all-cm-11}.
\section{Deduction of Theorem \ref{th:int-unrestricted} from Theorem~\ref{thma} }\label{se:remaining} 
Given an arbitrary cusp form, we first show how to construct an auxiliary form satisfying \ref{tag:1}, so that Theorem~\ref{thma} can be applied to it. The sign-change results of this auxiliary form then extend to the non-zero cusp form that we started with.

\subsection{Reduction to \texorpdfstring{\ref{tag:1}}{Hypothesis H}}\label{se:proof-D}
Let $0\neq f\in S_\kappa^*(\GammaOne{N})$  with Fourier expansion \eqref{expansion} such that $a_f(n)= 0$ whenever $(n,2N)=1.$ Then there is some $n_0> 1$ such that $a_f(n_0)\neq 0$. 
We set $N_f= \prod_{p|N} p$ and write $n_0= m_0 t$ with $(m_0, 2Nt)=1$. 
Such a factorisation always exists, one can simply take $m_0=1$, $t=n_0$, for example. 

Now for all primes $p\mid 2N_f$,  applying the prime-sifting operators $\mathcal{D}_p$ (defined in \eqref{eq:dp-operator}) to $f\mid_\kappa U_t$, we obtain
\[ g_t :=  f  \mid_\kappa U_t  \mid_\kappa \Big(\prod \nolimits_{p|2N_f}  \mathcal{D}_p \Big) ; \quad \text{ so that } \quad g_t(\tau) = \sumn_{m\ge 1, \, (m,2N_f)=1} a_f(mt) q^m.\] 
Note that $0\neq g_t\in S_\kappa (\Gamma_1(tNN^2_f))$ since $a_f(m_0 t) \neq 0.$ The orthogonality of $g_t$ and unary theta series follows from the following lemma.
In particular, $g_t \in S_\kappa^* (\Gamma_1(tNN^2_f))$ satisfies \ref{tag:1}, and Theorem \ref{thma} applies. Further, any lower bound for $\mathscr{S}_{g_t}^\flat(X)$ implies the same for $\mathscr{S}_f(X)$. \qed

\begin{lem}\label{le:duke-bound}
Suppose $f \in S_{3/2}^*(\Gamma_1(N))$, and let $h = \sum_{n\ge1} \alpha(n) a_f(tn) q^n$, with fixed $t$ and $\sup_n |\alpha(n)| \ll 1$. 
Assume that $h$ is a cusp form of weight $3/2$ of some level $L$. Then $h\in S_{3/2}^*(\Gamma_1(L))$.
\end{lem}
\begin{proof}
Duke's bound (see \cite{duke88}) for cusp forms in $S_{3/2}^*(\Gamma_1(N))$ gives $|a_h(n)| =  |\alpha(n) a_f(tn)|  \ll_{f,\varepsilon} n^{13/28+\varepsilon}$. 

We decompose the cusp form $h$ as $h = h_0 + h_1$ with $h_0\in \mathcal{U}_{3/2}(L)$ and $h_1\in S_{3/2}^*(\Gamma_1(L))$. Since Duke's bound holds for both $h$ and $h_1$, the same is true for $h_0$ too. However, if $h_0\neq 0$ then this is a contradiction, since a non-zero linear combination of unary theta series has coefficients of size $\asymp\sqrt n$ infinitely often along a suitable fixed square class. We refer to the proof of \cite[Lemma 2.1]{fun18} for details of the last argument.
\end{proof}

\subsection{Proof of Theorem \ref{th:int-unrestricted} for  \texorpdfstring{$\kappa\in \frac{1}{2}\Z$}{κ ∈ 12 Z}}
For $f\in S_{\kappa}^*(\Gamma_1(N))$ $\kappa\ge 3/2$, Theorem \ref{thma}, applied to $g_t$ on $[X/t,MX/t]$ for suitable $t$, yields the half-integral weight case of Theorem \ref{th:int-unrestricted}. 

However, the mixed case of weight $3/2$, namely $f=f_0+f_1$ with both $f_0\in \mathcal{U}_{3/2}(N)$ and $f_1\in S_{3/2}^*(\Gamma_1(N))$ are non-zero, needs a little care. From Section \ref{se:theta}, we have $\mathscr{N}_{f_0}(X) \asymp \sqrt{X}$. Consequently, $a_f(n) = a_{f_1}(n)$ for all but $O_f(\sqrt{X})$ many $n\in [X,MX]$. Therefore, 
\[\mathscr{S}_f(X;M)\ge \mathscr{S}_{f_1}(X;M)-O_f(\sqrt X)\gg X^{1-\eps}\] by Theorem \ref{thma}.
\qed

\subsection{Proof of Theorem \ref{th:int-unrestricted} for  \texorpdfstring{$\kappa\in \Z$}{κ ∈ Z}}
Recall that, for $\kappa\in\Z$, $f$ decomposes into a linear combination of dilates of primitive forms $f_j$. 
We recall that
\begin{equation}\label{eq:newform-decomp-simplified-2}
    f(\tau)=\sumn_{1\le j \le s} \sumn_{d\mid (N/N_j)}\alpha_{j,d}\,f_j(d\tau) \quad \text{with} \quad f^\circ (\tau) = \sumn_{1\le j\le s} \alpha_{j} ~f_j(\tau).
\end{equation}
Because of the coprimality condition in $ \mathscr S_f^{\flat}$(X), Theorem \ref{thma} depends only on the lacunarity of  $f^\circ$. 

\subsubsection{\texorpdfstring{$f$}{f} is non-lacunary}\label{se:non-lacunary-derivation}
First, suppose that $f$ is non-lacunary. Note that here we must have $\kappa \ge 2$.
Next, we consider the case when $f$ is non-lacunary, so there is some component $f_j$ in \eqref{eq:newform-decomp-simplified-2}  which is non-CM. 
Let $d_0$ be the least $d$ in \eqref{eq:newform-decomp-simplified-2} for which $\alpha_{j,d}\ne0$ for a non-CM $f_j$. 
If $d_0=1$, then $f^\circ$ is non-lacunary, and this case has already been treated in Section \ref{se:one-non-cm}. 
Hence, we may assume that $d_0>1$.

By the choice of $d_0$, there is an index $j_0$ such that $f_{j_0}$ is non-CM and $\alpha_{j_0,d_0}\neq 0$, while $\alpha_{j,d}=0$ for every non-CM $f_j$ and every $d<d_0$.
Then for a positive integer $m$ with $(m,N)=1$, we have 
\begin{align}
    a_f(d_0m) &={} \sumn_j \sumn_d \alpha_{j,d}\,a_j(d_0m/d) ={} \sumn_j \Big(\sumn_{d\mid d_0} \alpha_{j,d} ~ a_j \left( \tfrac{d_0}{d} \right)\Big) a_j(m) = \sumn_j \gamma_j a_j(m). \quad(\text{say})
\end{align}
Note that  by the minimality of $d_0$,  $\alpha_{j_0,d}=0$ if $d<d_0$. Therefore,
\begin{align}
    \gamma_{j_0} = \sumn_{d\mid d_0}\alpha_{j_0,d} ~ a_j\left(\tfrac{d_0}{d}\right) = \alpha_{j_0,d_0} \neq 0.
\end{align}
Now set $g:=f\mid_\kappa U_{d_0}\in S_\kappa(\Gamma_1(N))$, so that $a_g(n)=a_f(d_0n)$ for all $n$. 
If $g=0$, then we have $\sum_j\gamma_j a_j(p)=0$ for all primes $p\nmid 2N$. 
Since the $f_j$ are distinct newforms, strong multiplicity one implies that their corresponding unramified Hecke eigencharacters are pairwise distinct. But by Artin's lemma on the linear independence of distinct characters, applied to the (unramified) Hecke algebra, we must have $\gamma_{j_0} = 0$, a contradiction. (The reader may notice that one could also have used an orthogonality argument via the Rankin-Selberg $L$-functions.)
Thus $a_g(p)\neq 0$ for some prime $p\nmid 2N$. In particular, $g$ satisfies \ref{tag:1}. 
Moreover, the form $g^\circ$ involves the non-CM newform $f_{j_0}$ with coefficient $\gamma_{j_0}\neq 0$, so $g^\circ$ is non-lacunary.
Therefore, Theorem~\ref{thma} applies to $g$, and together with $a_g(m)=a_f(d_0m)$ it yields
\[\mathscr{S}_f(X)\gg X \asymp \mathscr{N}_f(X).\]

 \subsubsection{\texorpdfstring{$f$}{f} is lacunary} 
 Now suppose that $f$ is lacunary.

\noindent { $\bullet~ \kappa\ge 2$:}  If $f$ satisfies \ref{tag:1}, then $0 \neq f^\circ$ is lacunary, and from \eqref{eq:serre-estimates} and Section \ref{se:all-cm} we obtain the desired estimate for $\mathscr{S}_f(X)$.

If $f$ does not satisfy \ref{tag:1}, Section \ref{se:proof-D} provides a (lacunary) cusp form $g_t$ which does. 
 Now applying Theorem \ref{thma} to $g_t$ on the interval $[X/t,2X/t]$ gives
\begin{align}
      \mathscr S_f(X)\gg  X(\log X)^{-1/2} \asymp \mathscr{N}_f(X).
\end{align}

\noindent {$\bullet ~\kappa = 1$:} We proceed  by choosing the least $d_0$ such that $\min \{d: \alpha_{j,d}\neq 0, \theta_j = \theta_f\}$ and consider the cusp form $g= f\mid_1 U_{d_0}$. Then the same argument as in Section \ref{se:non-lacunary-derivation} gives \ref{tag:1} for $g$ with $\theta_{g^\circ} = \theta_f$ (see \eqref{eq:notation-wt1}). Thus, applying Theorem \ref{thma} to $g$, we have $\mathscr{S}_g^\flat(X) \gg_f X(\log X)^{-\theta_f}$. Replacing $X$ by $X/d_0$, we have the claimed estimate in  for $f$. This completes the proof of Theorem \ref{th:int-unrestricted}. \qed

We conclude this section with the following corollary, which follows immediately from the results of Sections \ref{se:half-int}, \ref{se:one-non-cm}, and \ref{se:all-cm-11}.
\begin{cor}\label{cor:nonvanishing}
    Let $N \ge 1$, $\kappa\in \frac{1}{2}\mathbb Z $ with $\kappa>\frac{3}{2}$, and $f \in S_\kappa (\Gamma_1(N))$  be as in \eqref{expansion} satisfying \ref{tag:1}. 
    Then for any $\eps>0$ and sufficiently large $X,$ we have 
    \begin{align}
        \mathscr N_f^{\flat}(X;M) \gg 
    \begin{cases}
    \begin{aligned}
        &X^{1-\eps}, &&  \text{ if } \kappa \in \tfrac12 \Z,\\ 
        &X(\log X)^{-1/2}, && \text{ if } \kappa \in \Z.
        \end{aligned}
\end{cases} 
\end{align} 
Moreover, for $\kappa$ integral, if $f^\circ$ is non-lacunary, then $ \mathscr N_f^{\flat}(X) \gg X.$
\end{cor}
 This improves the previously best-known lower bounds for $ \mathscr N_f^{\flat}(X)$, due to \cite{ana19, bd22}, viz., $N_f^{\flat}(X) \gg X^{1-\eps}$ in integral weight and $N_f^{\flat}(X) \gg X^{5/8-\eps}$ in half-integral weight case.

\section{Prime coefficients}\label{prime}
In this section, we prove Theorem \ref{th:prime-intro}, namely, we show that $ \displaystyle \mathscr{S}_f^\mathrm{pr} (X) \gg X^{1-\vartheta_0-\eps}.$
We prove the unconditional part of Theorem \ref{th:prime-intro}  in Section \ref{se:uncond-prime} and the conditional part under GRH in Section \ref{se:cond-prime}. 
To set-up the proofs, we first start with the following calculation. 
Note that, from \eqref{eq:newform-decomp-simplified-coeffs}, we have
\begin{align}
\sum_{\substack{x < p \le x+y\\ p \nmid N}} |c_{f}(p)|^2 \log p &= \sum_i |\alpha_i|^2 \sum_{\substack{x < p \le x+y\\ p \nmid N}} |\lambda_i(p)|^2 \log p + \sum_{i\neq j} \alpha_i \overline{\alpha_j} \sum_{\substack{x < p \le x+y\\ p \nmid N}} \lambda_i(p) \overline{\lambda_j(p)} \log p \nonumber\\
&= \sum_i |\alpha_i|^2 \sum_{\substack{x < p \le x+y\\ p \nmid N}} c_{f_i \times \Bar f_i}(p) + \sum_{i\neq j} \alpha_i \overline{\alpha_j} \sum_{\substack{x < p \le x+y\\ p \nmid N}} c_{f_i \times \Bar f_j}(p) , \label{linearcomb}
\end{align}
where $ c_{f_i \times \Bar f_j}(n) = \Lambda(n) \left( \sum\limits_{u,v} \alpha_{u,p;f_i}^m \, \Bar{\alpha}_{v,p;f_j}^m \right)$ whenever $n = p^m,$ and $0$ otherwise. 
Here $\alpha_{u,p;f}$ ($u=1,2$) denotes the $p$-Satake parameter of an eigenform $f$.
Then we have the logarithmic derivative of the Rankin-Selberg convolution of $f_i$ and $f_j$ as 
\begin{align}
    -\frac{L'}{L}(s,f_i\times\Bar{f}_j) = \sum_{n\geq 1} \frac{ c_{f_i \times \Bar f_j}(n)}{n^s}, \qquad c_{f_i \times \Bar f_j}(p) = \lambda_i(p) \overline{\lambda_j(p)} \log p.
\end{align}

\subsection{Unconditional result}\label{se:uncond-prime}
In this section we will prove the following variant of a theorem by Murty {\cite[Theorem 5]{Mur83}}.\footnote{There seems to be a mistake in \cite[Theorem 5]{Mur83} about the signs of Fourier coefficients; an arbitrary cusp form may not possess non-zero `prime' Fourier coefficients at all.}

\begin{thm}\label{unconditional:prime}
Let $f\in S_k(\Gamma_1(N))$ ($k \in \Z$) be a cusp form satisfying \ref{tag:1} with the Fourier expansion \eqref{expansion}. 
Then there exists a constant $0<\vartheta_0<1$ depending on $f$ such that for any $\eps>0$ and sufficiently large $X,$ we have $\displaystyle \mathscr S_f^{\mathrm{pr}}(X) \gg X^{1-\vartheta_0 - \eps}$.
\end{thm}

To prove Theorem \ref{unconditional:prime}, we need the following auxiliary lemma. 
Although Murty proves the following for level 1 eigenforms using \cite{mor73}, the same argument applied to the standard and Rankin-Selberg $L$-functions attached to eigenforms in $S_k(N,\chi)$ yields the following. 
\begin{lem}[{\cite[Lemma 4, 5(iii)]{Mur83}}]\label{murtylemma}
    Let $f,g\in  S_k(N,\chi)$ be distinct newforms with normalised Fourier coefficients $\lambda_f(n)$ and $\lambda_g(n)$, respectively. 
    Then there is a $\vartheta_0>0$ such that if $h=x^\vartheta,$ $\vartheta_0<\vartheta<1$, we have  
    \begin{enumerate}
        \item $\sum_{x \le p\le x+h} \lambda_f(p)\log p =o(h),$ \label{eq:1-1}
        \item $\sum_{x\le p\le x+h}\lambda_f(p)\overline{\lambda_g(p)}\log p=o(h),$ and \label{eq:1-2}
        \item $\sum_{x\le p\le x+h}|\lambda_f(p)|^2\log p\gg h.$ \label{eq:1-3}
    \end{enumerate}
\end{lem}

\begin{proof}[\textbf{Proof of Theorem \ref{unconditional:prime}}]
We first decompose $f$ into newforms $f_i$ of levels dividing $N$, as in \eqref{eq:newform-decomp-simplified}. 
Then applying \ref{eq:1-1} to each $f_i$  gives 
\begin{equation}\label{upp_uncond}
    \sumn_{x<p\le x+h} c_f(p)\log p=o_f (h),
\end{equation} 
for $h=x^\vartheta,$ with $\vartheta_0<\vartheta<1$. 
On the other hand, using \ref{eq:1-2}, \ref{eq:1-3} and the decomposition \eqref{linearcomb}, we obtain, for some constant $C_f>0$, $\sum_{\substack{x<p\le x+h \\ p\nmid N}} |c_f(p)|^2\log p \ge C_f h$ for all sufficiently large $x$. 
Since $c_f(p)\ll_f 1$, it follows that $\sum_{\substack{x<p\le x+h \\ p\nmid N}}~ |c_f(p)|\log p\gg_f h.$ 
Consequently, for all sufficiently large $x$,
\begin{equation}
\Big| \sumn_{\substack{x<p\le x+h,~ p\nmid N}}~c_f(p)\log p \Big| < \sumn_{\substack{x<p\le x+h,~ p\nmid N}} ~|c_f(p)|\log p,
\end{equation}
and the interval $(x,x+h]$ must contain primes $p\nmid N$ such that $c_f(p)$ takes both signs.  
Finally, partitioning $[X,2X]$ into disjoint intervals of length $\asymp X^{\vartheta_0+\eps}$ yields $\mathscr S_f^{\mathrm{pr}}(X) \gg X^{1-\vartheta_0 - \eps}.$
\end{proof}

\subsection{Under the assumption of GRH}\label{se:cond-prime} 
In this section, we assume the Generalised Riemann Hypothesis (GRH) for the automorphic $L$-functions associated with eigenforms and their Rankin-Selberg convolutions. 
Namely, all non-trivial zeros of these $L$-functions in the critical strip lie on the critical line $\Re(s)=\tfrac{1}{2}$.

\begin{thm}\label{signintegralprime}
    Assume the Generalised Riemann Hypothesis and let $f$ be as in Theorem \ref{unconditional:prime}. 
    Then 
    \[\mathscr S_f^{\mathrm{pr}}(X) \gg Xy^{-1}\] 
    for any increasing function $y=y(x) \le x$ with $y(\log x)^{-2} \rightarrow \infty$ as $x\rightarrow \infty.$

\end{thm}
 In particular, choosing $y=(\log x)^{2+\eps}$ gives the desired lower bound in Theorem \ref{th:prime-intro}.

\emph{Throughout the section, we fix the cusp form $f$ and $y=y(x)$ as in Theorem \ref{signintegralprime}.}

The proof of the theorem hinges on Propositions~\ref{prop-upperbound-prime} and \ref{prop-lowerbound-prime}, which are consequences of a common mean-square prime number theorem in short intervals. 
We therefore begin with the following lemma.
\begin{lem}\label{prop-primenumber}
    Assume the Generalised Riemann Hypothesis for $L(s,f_i)$ and $L(s,f_i\times\Bar{f}_j)$ for every $f_i,f_j$ in \eqref{eq:newform-decomp-simplified-coeffs}. Then,
    \begin{equation}\label{primenumberf}
    \int_X^{2X} \Big| \sum_{\substack{x< n\leq x+y,~ (n,2N)=1}} \Lambda(n) \lambda_i(n) \Big|^2 dx=o(Xy^2).
    \end{equation}
    and
    \begin{equation}\label{primenumberfg}
        \int_X^{2X} \Big| \sum_{\substack{x< n\leq x+y,~ (n,2N)=1}} c_{f_i \times \Bar f_j}(n)-\delta_{i,j}\,y\Big|^2 dx=o(Xy^2).
    \end{equation}
\end{lem}
\begin{proof}
    Our first assertion follows directly from \cite[Theorem 1.1]{Qu12}. 

    The Rankin-Selberg $L$-function $L(s,f_i\times \Bar f_j)$ satisfies the same analytic properties used in the proof of \cite[Theorem 1.1]{Qu12}, namely analytic continuation, a functional equation, the required zero-density estimates, and suitable bounds for its logarithmic derivative. 
    At every prime $p\nmid N$, Deligne's bound gives
    \[\Big| c_{f_i \times \Bar f_j}(p^m)\Big| =\Big|\log p \cdot \sumn_{1\le u,v\le 2} \alpha_{u,p;f_i}^{\,m} \bar{\alpha}_{v,p;f_j}^{\,m} \Big| \le 4\log p.\] 
    Hence the error in the Perron formula can be estimated directly from $|c_{f_i \times \Bar f_j}(n)|\ll_{f,g}\Lambda(n).$ 
    Shifting the Perron contour then  gives
    \[ \sumn_{n\le x}c_{f_i \times \Bar f_j}(n)=\delta_{i,j}x -\sumn_{|\gamma|\le T}\frac{x^\rho}{\rho}+\mathcal T_{i,j} (x,T) + O(R_{i,j}(x,T)),\]
    where $\rho = \beta +i\gamma$ runs over the non-trivial zeros of $L(s,f_i\times \Bar f_j)$, $\mathcal T_{i,j} (x,T)$ denotes the contribution of the trivial zeros, and $R_{i,j}$  is the error term.
    The only additional term compared to \cite{Qu12} is $\delta_{i,j} x$, arising from the simple pole at $s=1$ in the diagonal case. 
    Subtracting the formula at $x+y$ and $x$, and then following the Gallagher argument as in \cite[\S 4]{Qu12}, we obtain under GRH 
    \begin{align}
        \int_X^{2X}\Big|\sumn_{x< n\le x+y}c_{f_i \times \Bar f_j}(n) - \delta_{i,j}y \Big|^2 dx \ll y^2T\log^2(QT) + \frac{X^2\log^2(QT)}{T}+o(Xy^2),
    \end{align}
    where $Q$ is the conductor of $L(s,f_i \times\Bar{f_j})$. 
    Choosing $T=\frac{X}{y}$ gives 
    \begin{align}
        \int_X^{2X} \Big| \sumn_{x< n\le x+y} c_{f_i \times \Bar f_j}(n) - \delta_{i,j} y \Big|^2 dx \ll Xy\log^2 X + o(Xy^2) = o(Xy^2), 
    \end{align}
    since $y/(\log X)^2\to\infty$. 
\end{proof}

With these mean-square estimates in hand, Chebyshev's inequality allows us to pass from $L^2$-estimates to pointwise estimates for almost all $x\in[X,2X]$.
We first use this to prove the required cancellation in the signed prime sum.

\begin{prop}\label{prop-upperbound-prime}
      Assume GRH for every $f_i$ appearing in \eqref{eq:newform-decomp-simplified-coeffs}. Then for sufficiently large $X,$ the estimate 
\begin{equation}\label{upperboundsigned:prime1}
    \sumn_{\substack{ x< p\leq x+y,~ p\nmid N}} c_f(p) \log p  = o_f\left(y\right)
\end{equation} 
holds for almost all  $x\in [X,2X].$
\end{prop}
\begin{proof}
    For a cusp form $f\in S_\kappa(\Gamma_1(N))$, let
    \[\Delta(x,y;f):= \sumn_{\substack{x < n \le x + y,~  (n,2N)=1}} \Lambda(n) c_f(n),\] 
    where $y$ is as in Theorem \ref{signintegralprime}. 
    Then, for all $1\le j \le s$,  from \eqref{primenumberf} we have $\int_{X}^{2X} |\Delta(x,y;f_j)|^2 \,dx = o(y^2X).$ 
    Now applying Cauchy-Schwarz to $\Delta(x,y;f)=\sum_{j=1}^{s}\alpha_j\Delta(x,y;f_j)$ we have
    \begin{equation}\label{meansquare:lsfp}
        \int_{X}^{2X} |\Delta(x,y;f)|^2 \,dx = o(y^2X).
    \end{equation} 
Now, for any $\delta>0,$ we define  $E_\delta:= \Big\{ X\leq x\leq 2X ~:~  |\Delta(x,y;f)| \geq \delta y \Big\}.$ 
Then, from the definition of $E_\delta$, we have 
\begin{align*}
    \int_{X}^{2X} |\Delta(x, y; f)|^2 \, dx \geq (\delta y)^2 \cdot \mathrm{meas}(E_\delta).
\end{align*} 
To put it in other way, we have $ \mathrm{meas}(E_\delta) = o(X).$
Therefore $\Delta(x,y;f) = o(y) $ holds for almost all $x\in [X,2X].$ 
Since $\Lambda(n)$ is supported on prime powers, we can decompose $\Delta(x,y;f) $ as
\begin{equation}\label{mangoldtseparation:signed}
    \Delta(x,y;f) = \sum_{\substack{x< p\le x+y,~p\nmid N}} c_f(p) \log p + \sum_{\substack{x< p^r\le x+y\\r\ge2,~ p\nmid N}} c_f(p^r) \log p. 
\end{equation}

As $\Delta(x,y;f)$ and the second sum on the right-hand side are $o(y)$ for almost all $x\in [X,2X]$, we immediately have the assertion.
\end{proof}

We now turn to the complementary lower bound. Here the Rankin-Selberg estimate in \eqref{primenumberfg} is used to control the second moment of the prime coefficients: the diagonal terms provide the main term, while
the off-diagonal terms are negligible for almost all short intervals.

\begin{prop}\label{prop-lowerbound-prime}
     Assume GRH for $L(s,f_i\times\Bar{f}_j)$ for every pair $f_i,f_j$ appearing in \eqref{eq:newform-decomp-simplified-coeffs}. 
     Then 
    \begin{equation}\label{lowerboundsigned:prime}
        \sumn_{\substack{x<p\leq x+y,~ p\nmid N }}~ |c_f(p)| \log p \gg_f y
    \end{equation} 
    holds for almost all  $x\in [X,2X].$
\end{prop}
\begin{proof}
    We prove that $\sum_{\substack{x<  p\leq x+y\\ p\nmid N}} ~c_f(p)^2 \log p \gg_f y$ for almost all $x\in [X,2X].$ 
    Proceeding as in Proposition \ref{prop-upperbound-prime}, from \eqref{linearcomb} we conclude that for some $\alpha_i\neq 0$ and $C_i'>0$,
\begin{align}
    \sum_{\substack{x<p\leq x+y\\ p\nmid N }} |c_f(p)|^2 \log p &\geq |\alpha_i|^2 \sum_{\substack{x < p \le x+y\\ p \nmid N}} |\lambda_i(p)|^2 \log p + \sum_{i\neq j} \alpha_i \overline{\alpha_j} \sum_{\substack{x < p \le x+y\\ p \nmid N}} \lambda_i(p) \overline{\lambda_j(p)} \log p = C_i 'y + o(y)
\end{align} 
holds for almost all  $x\in [X,2X]$ and the claim follows. 
\end{proof}

\subsubsection{Proof of Theorem \ref{signintegralprime}}
    From Propositions \ref{prop-upperbound-prime} and  \ref{prop-lowerbound-prime}, we already have 
    \[\Big|\sumn_{\substack{x < p \leq x+y,~ p\nmid N}} c_f(p) \log p \Big| < \sumn_{\substack{x<p\leq x+y,~ p\nmid N }} |c_f(p)| \log p\] 
    for almost all $X\leq x\leq 2X$. 
    Consequently, Lemma \ref{le:sign-change-count} applies and we conclude that there are at least $\gg Xy^{-1}$ sign changes in $[X,2X].$ \qed

\section{Siegel cusp forms} \label{sec:smf}
We begin this section by recalling the basic notation and definitions for Siegel modular forms that will be used subsequently. The notation is standard, and further details can be found in \cite{kli90}. We denote the Siegel upper half‑space of degree $n$ by 
\[\mathbb H_n := \{ Z \in M_n(\mathbb C) \mid Z = Z^t,\ \operatorname{Im}(Z) > 0 \}.\]  
The symplectic group $\Sp(2n, \R) := \{g \in \GL(2n, \R) \; | \; J_n[g] = J_n \}$ acts on $\mathbb H_n$ by  
\[Z \mapsto g\langle Z \rangle = (AZ+B)(CZ+D)^{-1}, \qquad g = \Big(\begin{smallmatrix} A & B \\ C & D \end{smallmatrix}\Big) \in \mathrm{Sp}(2n, \mathbb R).\]  

A Siegel modular form of degree $n$ and weight $k \in \mathbb Z$ is a holomorphic function $F : \mathbb H_n \to \mathbb C$ satisfying  
\[F(g\langle Z \rangle) = \det(CZ+D)^k \, F(Z)\]
for all $g \in \mathrm{Sp}(2n, \mathbb Z)$, with the usual cusp condition when $n=1$.  
We denote by $M_k^n$ the space of Siegel modular forms of weight $k$ and degree $n.$ 
Every $F \in M_k^n$ admits a Fourier expansion  
\[F(Z) = \sumn_{S \in \Lambda_n} a_F(S) \, e(\mathrm{tr}(SZ)),\]
and an element $F\in M_k^n$ is called a cusp form if the sum is restricted to positive definite $S$, i.e., if $a_F(S) = 0$ unless $S\in \Lambda_n^+$. 
We denote by $S_k^n \subset M_k^n$ the subspace of all Siegel cusp forms.  

\subsection{Fourier-Jacobi expansion}
We now fix a decomposition of $Z \in \mathbb H_n$ into blocks  
\[Z = \begin{pmatrix} \tau & \z \\ \z^t & Z' \end{pmatrix}, \qquad \tau \in \mathbb H_1,\ \z \in \mathbb C^{1,n-1},\ Z' \in \mathbb H_{n-1}.\]
For $F \in S_k^n$, its Fourier-Jacobi expansion with respect to this decomposition is  
\[F(Z) = \sumn_{T \in \Lambda_{n-1}} \varphi_T(\tau, \z) \, e(\mathrm{tr}(TZ')),\]  
where each $\varphi_T$ is a Jacobi form of weight $k$ and index $T$. 
Details on Fourier-Jacobi expansions can be found in \cite[\S 8]{kli90}, and for Jacobi forms of matrix index we refer to \cite{zie89}.

A Jacobi form $\varphi$ of weight $k$ and index $T$ admits a theta expansion  
\[\varphi(\tau, \z) = \sumn_{\mu \in \mathbb Z^{n-1} / 2T \mathbb Z^{n-1}} h_\mu(\tau) \, \Theta_T[\mu](\tau, \z),
\]  
where the functions $h_\mu$ are elliptic modular forms of weight $k - \frac{n-1}{2}$ on some congruence subgroup (see \cite{zie89}). 
The Fourier expansion of $h_\mu$ is given by  (see \cite[\S 2.3]{bd22})
\begin{equation}\label{hmu}
    h_\mu(\tau) = \sumn_{m \ge 0} a_\mu\!\left(m - T^{-1}[\mu/2]\right) q^{m - T^{-1}[\mu/2]}.
\end{equation} 

In particular, when $\varphi_T$ is the $T$-th Fourier-Jacobi coefficient of $F\in S_k^n$, then the Fourier coefficients of $h_\mu$ are related to those of $F$ by
\begin{equation}\label{eqcoeff}
    a_\mu\!\left(m - T^{-1}[\mu/2]\right) = a_F\!\left( \begin{pmatrix} m & \mu/2 \\ \mu^t/2 & T \end{pmatrix} \right).
\end{equation} 

Recall that for $M \in \Lambda_m^+$, a $\mu\in \mathbb Z^m$ is called primitive if $\frac{1}{4} M^{-1}[\mu] $ has exact denominator $d_M$ if $m$ is even and $4d_M$ if $m$ is odd. 
By  \cite[Proposition 2.3]{bd22}, a primitive $\mu$ exists whenever $d_M$ is odd and square-free.

\subsection{Proof of Theorem~\ref{siegelpositiveprop} -- reduction to one variable}
For $F \in S_k^n$ non-zero, we first obtain (infinitely many) $T \in \Lambda_{n-1}^+$ with $d_T$ odd and square-free such that $\varphi_T \neq 0$. This is worked out in \cite{bd22}. One can simply apply \cite[Theorem~1.1]{bd22} to the degree $n-1$ (vector-valued) form $F^\circ$ obtained through a Taylor expansion of $F$, we do not recall it here as it will not be used any further. See \cite[Section~3, Proposition~3.1]{bd22} for the construction of $F^\circ$. Its Fourier coefficients are proportional to (a linear combination of) the derivatives  $\frac{\partial}{\partial z^\lambda}\phi_T(\tau, z)\vert_{z=0}$ , where $\lambda \in \Z^{n-1}$ are certain vectors. Thus $F^\circ \neq 0$ implies the existence of one odd, square-free $T \in \Lambda_{n-1}$ such that $\phi_T \neq 0$.

Next, by \cite[Proposition 3.5]{bd22}, the Fourier–Jacobi coefficient $\varphi_T$ has a non-zero theta component $h_\mu$ with $\mu$ primitive. 
This $h_\mu$ is an elliptic modular form of weight $\kappa = k-\frac{n-1}{2}$ on $\Gamma(d_T')$, where $d_T' = d_T$ for $n$ odd and $d_T' = 4d_T$ for $n$ even. 
Set \[f(\tau) = h_\mu(d_T'\tau) \in M_\kappa\bigl(\Gamma_1(d_T'^2)\bigr).\]
Since $\mu$ is primitive, $T^{-1}[\mu/2]$ has exact denominator $d_T'$, and hence so does $m - T^{-1}[\mu/2]$ for every integer $m \ge 1$. 
By \eqref{hmu} and \cite[Proposition~3.5]{bd22}, the Fourier expansion of $f$ is therefore supported away from its level:
\begin{equation}\label{away}
f(\tau) = \sum\nolimits_{\substack{m\ge 1,~ (m,d_T')=1}} a_f\left(m\right) q^m.
\end{equation} 

From \eqref{hmu}, \eqref{eqcoeff} and \eqref{away}, we obtain that $f$ satisfies \ref{tag:1} with Fourier coefficients given by 
\begin{equation}\label{eltosi}
    a_f(m) = a_F\left( \mathcal T_m \right), \quad \mathcal T_m = \begin{pmatrix}
        \frac{m}{d_T'} + T^{-1}[\mu/2] & \mu/2\\\mu^t/2 & T
    \end{pmatrix}
\end{equation}
where we make the convention that $a_F \left(\left( \begin{smallmatrix}
        a & *\\ * & *
    \end{smallmatrix} \right)\right) = 0$ if $a$ is not an integer. 
    Then it is easy to verify that $d_{\mathcal T_m} =m$ (see, for example, \cite[Section 3.4]{bd22}). Now we prove our main theorem for odd and square-free discriminants.

\subsection{Proof of Theorem \ref{siegelpositiveprop} -- working with one variable}
    By the above subsection, we continue to work with  $0\neq f\in S_\kappa(\Gamma_1(d_T'^2))$ satisfying \ref{tag:1}, defined in \eqref{away}.
    Note that $\kappa $ is always $ \ge 2$ if $\kappa \in Z$ and it is $\ge 5/2$ if $\kappa$ is half-integral under the assumptions we have.
    
    Since $F$ has real Fourier coefficients, the same holds for $f$. Invoking Theorem~\ref{thma}, we obtain
    \begin{align}\label{eq:siegel-combined}
        \mathscr{S}^{\flat}_f(X;M) \geq C_{f} X L(X),
    \end{align}
    for some positive constant $C_f > 0$. Here $L(X)$ is either $(\log X)^{-1/2}$ or $X^{-\eps}$, according as $\kappa$ is integral or half-integral. 
    

    Let $\kappa\in \frac 12 \Z$ and $X_0\ge 1$ such that $X_0^{\eps}\ge C_f^{-2}$. Then  for all $X\geq X_0,$  we have     $\mathscr{S}^{\flat}_f(X;M) \geq X^{1-\eps}$ (cf. Theorem \ref{thma}).
    
    On the other hand, for $\kappa$ integral, let $A>\max\{1,C_f^{-1}\}$. Now applying \eqref{eq:siegel-combined} with $AX$ in place of $X$, we obtain $\mathscr{S}^{\flat}_f(X;2A) \geq \mathscr{S}^{\flat}_f(AX) \geq X(\log X)^{-1/2}$ for sufficiently large $X$ (cf. Theorem \ref{thma}).

    Consequently, in either case, for sufficiently large $X$, there exist an integer $r \geq XL(X)$ and an increasing sequence $(m_i)_{1 \leq i \leq r}$ of odd and square-free integers such that $c_f(m_i) c_f(m_{i+1}) < 0$.  For each such $m_i$, from \eqref{eltosi} we see that $a_f(m_i) = a_F(\mathcal{T}_i)$ with $d_{\mathcal{T}_i}=m_i$. Thus, for sufficiently large $X$, the estimates for $\mathscr{S}(X,\flat)$ follow.
      
     When $n$ is odd, we apply  Theorem~\ref{unconditional:prime} and \ref{signintegralprime} to the cusp form $f$ in \eqref{away}. We use the 
     relation \eqref{eltosi} to transfer these sign changes to the Fourier coefficients of $F$ with prime discriminants, we obtain
     \begin{align}
         \mathscr{S}(X,\mathbb P) \ge \mathscr{S}_f^\mathrm{pr}(X) \gg X^{1-\delta-\eps},
     \end{align}
     for any $\eps>0$. 
     In fact, if we assume GRH, then we obtain
     \begin{align}
          \mathscr{S}(X,\mathbb P) \gg X(\log X)^{-(2+\eps)} \gg X^{1-\eps}
     \end{align}
     and the proof is complete.
     \qed

{\small \subsection*{Acknowledgements} 
The first author thanks IISc. Bangalore, DST (India) and UGC centre for advanced studies for financial support. 
The second named author gratefully acknowledges the financial support provided by the National Board for Higher Mathematics (NBHM), Department of Atomic Energy (DAE), Government of India, under the Post-Doctoral Fellowship scheme. 
S.D. is supported by an ANRF grant ANRF/ARGM/2025/000535/MTR from the Ministry of Science and Technology, India. 
This research is supported by the DST FIST program-2021 [TPN-700661].}

{\small \subsection*{AI usage declaration} 
    	During the preparation of this work, the authors used ChatGPT to assist with some computations in Sections \ref{sec:32alt} and \ref{se:upper-bound-signed}. The original ideas presented in this paper are those of the authors. Following the use of ChatGPT, the authors reviewed and edited the content as needed and take full responsibility for the results of the paper.
}
\printbibliography
\end{document}